\documentclass[12pt]{article}
\usepackage{amsfonts}\usepackage{amssymb}\usepackage{amsmath}

\newcommand{\ignore}[1]{}

\def\@begintheorem#1#2{\par\bgroup{\sc #1\ #2. }\it\ignorespaces}
\def\@opargbegintheorem#1#2#3{\par\bgroup{\sc #1\ #2\ (#3). } \it\ignorespaces}
\def\@endtheorem{\egroup}
\newtheorem{theorem}{Theorem}[section]
\newtheorem{corollary}[theorem]{Corollary}
\newtheorem{lemma}[theorem]{Lemma}
\newtheorem{proposition}[theorem]{Proposition}
\newtheorem{problem}[theorem]{Problem}
\newtheorem{example}[theorem]{Example}
\newtheorem{algorit}[theorem]{Algorithm}
\newtheorem{definition}[theorem]{Definition}
\newcommand{\bt}[1]{\begin{theorem}\label{#1}}
\newcommand{\bc}[1]{\begin{corollary}\label{#1}}
\newcommand{\bl}[1]{\begin{lemma}\label{#1}}
\newcommand{\bp}[1]{\begin{proposition}\label{#1}}
\newcommand{\bpro}[1]{\begin{problem}\label{#1}}
\newcommand{\be}[1]{\begin{example}\rm\label{#1}}
\newcommand{\ba}[1]{\begin{algorit}\rm\label{#1}}
\newcommand{\bd}[1]{\begin{definition}\rm\label{#1}}
\newcommand{\bpr}{\noindent {\em Proof. }}
\newcommand{\et}{\end{theorem}}
\newcommand{\ec}{\end{corollary}}
\newcommand{\el}{\end{lemma}}
\newcommand{\ep}{\end{proposition}}
\newcommand{\epro}{\end{problem}}
\newcommand{\ee}{\end{example}}
\newcommand{\ea}{\end{algorit}}
\newcommand{\ed}{\end{definition}}
\newcommand{\epr}{{\ \vbox{\hrule\hbox{%
\vrule height1.3ex\hskip0.8ex\vrule}\hrule
}}\\\par}

\def\N{\mathbb{N}}
\def\R{\mathbb{R}}
\def\Z{\mathbb{Z}}
\def \F {{{\cal F}}}
\def \G {{{\cal G}}}
\def \L {{{\cal L}}}
\def \l {\langle}
\def \r {\rangle}
\def \lin  {{\rm lin}}
\def \aff  {{\rm aff}}
\def \conv {{\rm conv}}
\def \cone {{\rm cone}}
\def \vert {{\rm vert}}
\def \supp {{\rm supp}}
\def \zone {{\rm zone}}
\def \sign {{\rm sign}}

\begin{document}

\thispagestyle{empty}
\begin{center}
\vskip2cm
{\bf\Huge Convex Discrete Optimization}
\vskip6cm
{\bf\LARGE Shmuel Onn}
\vskip1cm
{\large Technion - Israel Institute of Technology, Haifa, Israel}
\vskip9cm
{\large March 2007}
\end{center}

\newpage\ \vskip5cm
\thispagestyle{empty}

\centerline{\bf Abstract}
\vskip.5cm
We develop an algorithmic theory of convex optimization over discrete sets.
Using\break a combination of algebraic and geometric tools we are able to
provide polynomial time algorithms for solving broad classes of convex
combinatorial optimization problems and convex integer programming problems in
variable dimension. We discuss some of the many applications of this theory
including to quadratic programming, matroids, bin packing and cutting-stock
problems, vector partitioning and clustering, multiway transportation problems,
and privacy and confidential statistical data disclosure. Highlights of our
work include a strongly polynomial time algorithm for convex and linear
combinatorial optimization over any family presented by a membership
oracle when the underlying polytope has few edge-directions; a new theory of
so-termed $n$-fold integer programming, yielding polynomial time solution of
important and natural classes of convex and linear integer programming
problems in variable dimension; and a complete complexity classification of
high dimensional transportation problems, with practical applications to
fundamental problems in privacy and confidential statistical data disclosure.

\newpage
{\small \tableofcontents}

\newpage
\section{Introduction}
\label{introduction}

The general linear discrete optimization problem can be posed as follows.

\vskip.2cm\noindent
{\sc Linear discrete optimization}. Given a set $S\subseteq\Z^n$
of integer points and an integer vector $w\in\Z^n$, find an
$x\in S$ maximizing the standard inner product $wx:=\sum_{i=1}^nw_ix_i$.

\vskip.2cm\noindent
The algorithmic complexity of this problem, which includes
{\em integer programming} and
{\em combinatorial optimization} as special cases,
depends on the presentation
of the set $S$ of feasible points. In integer programming,
this set is presented as the set of integer points
satisfying a given system of linear inequalities,
which in standard form is given by
$$S\quad=\quad\{x\in\N^n:\ Ax=b\}\quad,$$
where $\N$ stands for the nonnegative integers, $A\in\Z^{m\times n}$
is an $m\times n$ integer matrix, and $b\in\Z^m$ is an integer vector.
The input for the problem then consists of $A,b,w$.\break
In combinatorial optimization, $S\subseteq\{0,1\}^n$ is a set of
$\{0,1\}$-vectors, often interpreted as a family of subsets of a ground
set $N:=\{1,\dots,n\}$, where each $x\in S$ is the indicator of its
support $\supp(x)\subseteq N$. The set $S$ is presented implicitly
and compactly, say as the set of indicators of subsets of edges
in a graph $G$ satisfying a given combinatorial property
(such as being a matching, a forest, and so on), in which case
the input is $G,w$. Alternatively, $S$ is given by an oracle, such
as a {\em membership oracle} which, queried on $x\in\{0,1\}^n$, asserts
whether or not $x\in S$, in which case the algorithmic complexity also
includes a count of the number of oracle queries needed to solve the problem.

Here we study the following broad generalization
of linear discrete optimization.

\vskip.2cm\noindent
{\sc Convex discrete optimization}.
Given a set $S\subseteq\Z^n$, vectors $w_1,\dots,w_d\in\Z^n$,
and a convex functional $c:\R^d\longrightarrow\R$, find an $x\in S$
maximizing $c(w_1 x,\dots,w_d x)$.

\vskip.2cm\noindent
This problem can be interpreted as {\em multi-objective}
linear discrete optimization: given $d$ linear functionals
$w_1x,\dots,w_dx$ representing the values of points $x\in S$
under $d$ criteria, the goal is to maximize their ``convex balancing"
defined by $c(w_1 x,\dots,w_d x)$. In fact, we have a hierarchy of problems
of increasing generality and complexity, parameterized by the number $d$
of linear functionals: at the bottom lies the linear discrete optimization
problem, recovered as the special case of $d=1$ and $c$ the identity
on $\R$; and at the top lies the problem of maximizing an arbitrary
convex functional over the feasible set $S$, arising with $d=n$ and with
$w_i={\bf 1}_i$ the $i$-th standard unit vector in $\R^n$ for all $i$.

The algorithmic complexity of the convex discrete optimization problem
depends on the presentation of the set $S$ of feasible points as in
the linear case, as well as on the presentation of the convex functional $c$.
When $S$ is presented as the set of integer points satisfying a
given system of linear inequalities we also refer to the problem as
{\em convex integer programming}, and when $S\subseteq\{0,1\}^n$
and is presented implicitly or by an oracle we also refer to
the problem as {\em convex combinatorial optimization}. As for the
convex functional $c$, we will assume throughout that it is presented
by a {\em comparison oracle} that, queried on $x,y\in\R^d$, asserts whether
or not $c(x)\leq c(y)$. This is a very broad presentation that reveals
little information on the function, making the problem, on the one hand,
very expressive and applicable, but on the other hand, very hard to solve.

There is a massive body of knowledge on the complexity of linear
discrete optimization - in particular (linear) integer programming
\cite{Sch} and (linear) combinatorial optimization \cite{GLS}.
The purpose of this monograph is to provide the first comprehensive
unified treatment of the extended convex discrete optimization problem.
The monograph follows the outline of five lectures given by the author
in the S\'eminaire de Math\'ematiques Sup\'erieures Series,
Universit\'e de Montr\'eal, during June 2006. Colorful slides of theses
lectures are available online at \cite{Onn} and can be used as a visual
supplement to this monograph. The monograph has been written under the
support of the ISF - Israel Science Foundation. The theory developed here
is based on and is a culmination of several recent papers including
\cite{BOT,DO1,DO2,DO3,DO4,DHOW,DHORW,FOR,HOR,Onn1,Onn2,OR,OSc,OSt}
written in collaboration with several colleagues - Eric Babson, Jesus De Loera,
Komei Fukuda, Raymond Hemmecke, Frank Hwang, Vera Rosta, Uriel Rothblum,
Leonard Schulman, Bernd Sturmfels, Rekha Thomas, and Robert Weismantel.
By developing and using a combination of geometric and algebraic tools,
we are able to provide polynomial time algorithms for several broad classes
of convex discrete optimization problems. We also discuss in detail some
of the many applications of our theory, including to quadratic programming,
matroids, bin packing and cutting-stock problems, vector partitioning
and clustering, multiway transportation problems, and privacy and
confidential statistical data disclosure.

We hope that this monograph will, on the one hand, allow users of discrete
optimization to enjoy the new powerful modelling and expressive capability
of convex discrete optimization along with its broad polynomial time
solvability, and on the other hand, stimulate more research on this new
and fascinating class of problems, their complexity, and the study of
various relaxations, bounds, and approximations for such problems.

\subsection{Limitations}
\label{l}

Convex discrete optimization is generally intractable even for small
fixed $d$, since already for $d=1$ it includes linear integer programming
which is NP-hard. When $d$ is a variable part of the input, even very
simple special cases are NP-hard, such as the following problem,
so-called {\em positive semi-definite quadratic binary programming},
$$\max\,\{(w_1 x)^2+\cdots+(w_n x)^2\ :
\ x\in\N^n\,,\ x_i\leq 1\,,\ i=1,\dots,n\}\ .$$
Therefore, throughout this monograph we will
assume that $d$ is fixed (but arbitrary).

As explained above, we also assume throughout that the convex functional $c$
which constitutes part of the data for the convex discrete optimization problem is
presented by a comparison oracle. Under such broad presentation,
the problem is generally very hard. In particular, if the feasible
set is $S:=\{x\in\N^n:\ Ax=b\}$ and the underlying polyhedron
$P:=\{x\in\R^n_+\,:\, Ax=b\}$ is {\em unbounded}, then the problem
is inaccessible even in one variable with no equation constraints.
Indeed, consider the following family of univariate convex integer programs
with convex functions parameterized by $-\infty< u\leq\infty$,
$$
\max\, \{c_u(x)\ :\ x\in\N\}\ ,\quad
c_u(x):=\left\{
\begin{array}{ll}
    -x, & \hbox{if   $x<u$;} \\
    x-2u, & \hbox{if $x\geq u$.} \\
\end{array}
\right.\ .
$$
Consider any algorithm attempting to solve the problem and let $u$ be the
maximum value of $x$ in all queries to the oracle of $c$. Then the algorithm
can not distinguish between the problem with $c_u$, whose objective function is
unbounded, and the problem with $c_\infty$, whose optimal objective value is $0$.
Thus, convex discrete optimization (with an oracle presented functional)
over an {\em infinite} set $S\subset\Z^n$ is quite hopeless.
Therefore, an algorithm that solves the convex discrete optimization
problem will either return an optimal solution, or assert that the
problem is infeasible, or assert that the underlying polyhedron is unbounded.
In fact, in most applications, such as in combinatorial optimization
with $S\subseteq \{0,1\}^n$ or integer programming with
$S:=\{x\in\Z^n:Ax=b,\ l\leq x\leq u\}$ and $l,u\in\Z^n$, the set $S$ is
finite and the problem of unboundedness does not arise.

\subsection{Outline and Overview of Main Results and Applications}
\label{oaoomraa}

We now outline the structure of this monograph and provide a brief overview
of what we consider to be our main results and main applications. The
precise relevant definitions and statements of the theorems and corollaries
mentioned here are provided in the relevant sections in the monograph body.
As mentioned above, most of these results are adaptations or extensions
of results from one of the papers
\cite{BOT,DO1,DO2,DO3,DO4,DHOW,DHORW,FOR,HOR,Onn1,Onn2,OR,OSc,OSt}.
The monograph gives many more applications and results that may turn out to
be useful in future development of the theory of convex discrete optimization.

The rest of the monograph consists of five sections.
While the results evolve from one section to the next, it is quite easy
to read the sections independently of each other (while just browsing now
and then for relevant definitions and results). Specifically,\break
Section \ref{ccoam} uses definitions and the main result of
Section \ref{rctldo}; Section \ref{cip} uses definitions and results from
Sections \ref{rctldo} and \ref{lnfip}; and Section \ref{mtpapisd} uses
the main results of Sections \ref{lnfip} and \ref{cip}.

\vskip.2cm
In Section \ref{rctldo} we show how to reduce the convex discrete
optimization problem over $S\subset\Z^n$ to strongly polynomially many
linear discrete optimization counterparts over $S$, provided that the
convex hull $\conv(S)$ satisfies a suitable geometric condition, as follows.

\vskip.2cm\noindent{\bf Theorem \ref{ConvexToLinear}}
{\em For every fixed $d$, the convex discrete optimization problem over
any finite $S\subset\Z^n$ presented by a linear discrete optimization oracle
and endowed with a set covering all edge-directions of $\conv(S)$,
can be solved in strongly polynomial time.}

\vskip.2cm\noindent
This result will be incorporated in the polynomial time algorithms
for convex combinatorial optimization and convex integer
programming to be developed in \S \ref{ccoam} and \S \ref{cip}.

\vskip.2cm
In Section \ref{ccoam} we discuss convex combinatorial optimization.
The main result is that convex combinatorial optimization over a set
$S\subseteq\{0,1\}^n$ presented by a membership oracle
can be solved in strongly polynomial time provided it is endowed
with a set covering all edge-directions of $\conv(S)$.
In particular, the standard linear combinatorial optimization
problem over $S$ can be solved in strongly polynomial time as well.

\vskip.2cm\noindent{\bf Theorem \ref{CCO}}
{\em For every fixed $d$, the convex combinatorial optimization problem over
any $S\subseteq\{0,1\}^n$ presented by a membership oracle and endowed with
a set covering all edge-directions of the polytope $\conv(S)$,
can be solved in strongly polynomial time.}

\vskip.2cm\noindent
An important application of Theorem \ref{CCO}
concerns convex matroid optimization.

\vskip.2cm\noindent{\bf Corollary \ref{Matroids}}
{\em For every fixed $d$, convex combinatorial optimization
over the family of bases of a matroid presented by membership
oracle is strongly polynomial time solvable.}

\vskip.2cm
In Section \ref{lnfip} we develop the theory of linear {\em $n$-fold integer
programming}. As a consequence of this theory we are able to solve a broad
class of linear integer programming problems in variable dimension in
polynomial time, in contrast with the general intractability of linear
integer programming. The main theorem here may seem a bit technical at a
first glance, but is really very natural and has many applications
discussed in detail in \S \ref{lnfip}, \S \ref{cip} and \S \ref{mtpapisd}.
To state it we need a definition.
Given an $(r+s)\times t$ matrix $A$, let $A_1$ be its $r\times t$ sub-matrix
consisting of the first $r$ rows and let $A_2$ be its $s\times t$
sub-matrix consisting of the last $s$ rows. We refer to $A$ explicitly
as {\em $(r+s)\times t$ matrix}, since the definition below depends also
on $r$ and $s$ and not only on the entries of $A$. The {\em $n$-fold matrix}
of an $(r+s)\times t$ matrix $A$ is then defined to be
the following $(r+ns)\times nt$ matrix,
$$A^{(n)}\quad:=\quad ({\bf 1}_n\otimes A_1)\oplus(I_n \otimes A_2)\quad=\quad
\left(
\begin{array}{ccccc}
  A_1    & A_1    & A_1    & \cdots & A_1    \\
  A_2  & 0      & 0      & \cdots & 0      \\
  0  & A_2      & 0      & \cdots & 0      \\
  \vdots & \vdots & \vdots & \ddots & \vdots \\
  0  & 0      & 0      & \cdots & A_2      \\
\end{array}
\right)\quad .
$$
Given now any $n\in\N$, lower and upper bounds
$l,u\in\Z_{\infty}^{nt}$ with $\Z_{\infty}:=\Z\uplus\{\pm\infty\}$,
right-hand side $b\in\Z^{r+ns}$, and linear functional $wx$ with
$w\in\Z^{nt}$, the corresponding linear $n$-fold integer programming
problem is the following program in variable dimension $nt$,
$$\max\,\{wx\ :\ x\in\Z^{nt},\ A^{(n)}x=b,\ l\leq x\leq u\}\ .$$
The main theorem of \S \ref{lnfip} asserts that such
integer programs are polynomial time solvable.

\vskip.2cm\noindent{\bf Theorem \ref{NFoldTheorem}\ }
{\em For every fixed $(r+s)\times t$ integer matrix $A$,
the linear $n$-fold integer programming problem with any
$n$, $l$, $u$, $b$, and $w$ can be solved in polynomial time.}

\vskip.2cm\noindent
Theorem \ref{NFoldTheorem} has very important applications to high-dimensional
transportation problems which are discussed in \S \ref{twlstp} and in more
detail in \S \ref{mtpapisd}. Another major application concerns
bin packing problems, where items of several types are to be packed
into bins so as to maximize packing utility subject to weight constraints.
This includes as a special case the classical cutting-stock
problem of \cite{GG}. These are discussed in detail in \S \ref{ppacs}.

\vskip.2cm\noindent{\bf Corollary \ref{Packing}}
{\em For every fixed number $t$ of types and type weights
$v_1,\dots,v_t$, the corresponding integer bin packing and
cutting-stock problems are polynomial time solvable.}

\vskip.2cm
In Section \ref{cip} we discuss convex integer programming, where the
feasible set $S$ is presented as the set of integer points
satisfying a given system of linear inequalities.\break
In particular, we consider convex integer programming over $n$-fold systems for
any fixed (but arbitrary) $(r+s)\times t$ matrix $A$, where, given $n\in\N$,
vectors $l,u\in\Z_{\infty}^{nt}$, $b\in\Z^{r+ns}$ and $w_1,\dots,w_d\in\Z^{nt}$,
and convex functional $c:\R^d\longrightarrow\R$, the problem is
$$\max\,\{c(w_1x,\dots,w_dx)\ :\ x\in\Z^{nt},\ A^{(n)}x=b,\ l\leq x\leq u\}\ .$$
The main theorem of \S \ref{cip} is the following extension of
Theorem \ref{NFoldTheorem}, asserting that convex integer programming
over $n$-fold systems is polynomial time solvable as well.

\vskip.2cm\noindent{\bf Theorem \ref{NFoldConvex}\ }
{\em For every fixed $d$ and $(r+s)\times t$ integer matrix $A$,
convex $n$-fold integer programming with any $n$, $l$, $u$, $b$,
$w_1,\dots,w_d$, and $c$ can be solved in polynomial time.}

\vskip.2cm\noindent
Theorem \ref{NFoldConvex} broadly extends the class of objective functions
that can be efficiently maximized over $n$-fold systems.
Thus, all applications discussed in \S \ref{salip} automatically extend
accordingly. These include convex high-dimensional
transportation problems and convex bin packing and cutting-stock problems,
which are discussed in detail in \S \ref{tpapp} and \S \ref{mtpapisd}.

Another important application of Theorem \ref{NFoldConvex} concerns vector
partitioning problems which have applications in many areas including load
balancing, circuit layout, ranking, cluster analysis, inventory, and
reliability, see e.g. \cite{BHR,BH,FOR,HOR,OSc} and the references therein.
The problem is to partition $n$ items among $p$ players so as to maximize
social utility. With each item is associated a $k$-dimensional vector
representing its utility under $k$ criteria. The social utility of a
partition is a convex function of the sums of vectors of items
that each player receives. In the constrained version of the problem,
there are also restrictions on the number of items each player
can receive. We have the following consequence of Theorem \ref{NFoldConvex};
more details on this application are in \S \ref{vpac}.

\vskip.2cm\noindent{\bf Corollary \ref{Partition}}
{\em For every fixed number $p$ of players and number $k$
of criteria, the constrained and unconstrained vector partition
problems with any item vectors, convex utility, and constraints
on the number of item per player, are polynomial time solvable.}

\vskip.2cm
In the last Section \ref{mtpapisd} we discuss multiway (high-dimensional)
transportation problems and secure statistical data disclosure. Multiway
transportation problems form a very important class of discrete optimization
problems and have been used and studied extensively in the operations research
and mathematical programming literature, as well as in the statistics
literature in the context of secure statistical data disclosure and management
by public agencies, see e.g. \cite{AT,BR,Cox,DT,DLTZ,KW,KLS,QS,Vla,YKK} and
the references therein. The feasible points in a transportation problem are
the multiway tables (``contingency tables" in statistics) such that the
sums of entries over some of their lower dimensional sub-tables such as lines
or planes (``margins" in statistics) are specified. We completely settle the
algorithmic complexity of treating multiway tables and discuss the applications
to transportation problems and secure statistical data disclosure, as follows.

In \S \ref{tut} we show that ``short" $3$-way transportation problems,
over $r\times c\times 3$ tables with variable number $r$ of rows and variable
number $c$ of columns but fixed small number $3$ of layers (hence ``short"),
are {\em universal} in that {\em every} integer programming problem
is such a problem (see \S \ref{tut} for the precise
stronger statement and for more details).

\vskip.2cm\noindent{\bf Theorem \ref{Universality}}
{\em Every linear integer programming problem $\max\{cy:y\in\N^n:Ay=b\}$ is
polynomial time representable as a short $3$-way line-sum transportation problem
$$\max\,\{\,wx\ :\ x\in\N^{r\times c\times 3}\ :\ \sum_i x_{i,j,k}=z_{j,k}
\,,\ \sum_j x_{i,j,k}=v_{i,k}\,,\ \sum_k x_{i,j,k}=u_{i,j}\,\}\ .$$}

In \S \ref{tcotmtp} we discuss $k$-way transportation problems of any
dimension $k$. We provide the first polynomial time algorithm
for convex and linear ``long" $(k+1)$-way transportation problems,
over $m_1\times\cdots\times m_k\times n$ tables, with $k$ and
$m_1,\ldots, m_k$ fixed (but arbitrary), and variable number $n$ of
layers (hence ``long"). This is best possible in view of
Theorem \ref{Universality}. Our algorithm works for any
{\em hierarchical collection of margins}: this captures common margin
collections such as all line-sums, all plane-sums, and more generally all $h$-flat
sums for any $0\leq h\leq k$ (see \S \ref{tam} for more details).
We point out that even for the very
special case of linear integer transportation over $3\times 3\times n$
tables with specified line-sums, our polynomial time algorithm is the
only one known. We prove the following statement.

\vskip.2cm\noindent{\bf Corollary \ref{ConvexKWay}}
{\em For every fixed $d,k,m_1,\dots,m_k$ and family $\F$ of subsets of
$\{1,\dots,k+1\}$ specifying a hierarchical collection of margins,
the convex (and in particular linear) long transportation problem
over $m_1\times\cdots\times m_k\times n$ tables is polynomial time solvable.}

\vskip.2cm
In our last subsection \S \ref{paeu} we discuss an important
application concerning privacy in statistical databases.
It is a common practice in the disclosure of a multiway table containing
sensitive data to release some table margins rather than the
table itself. Once the margins are released, the security of any specific
entry of the table is related to the set of possible values that
can occur in that entry in any table having the same margins as those
of the source table in the data base. In particular, if this set consists of
a unique value, that of the source table, then this entry can be exposed and
security can be violated. We show that for multiway tables where one category
is significantly richer than the others, that is, when each sample point can
take many values in one category and only few values in the other categories,
it is possible to check entry-uniqueness in polynomial time, allowing
disclosing agencies to make learned decisions on secure disclosure.

\vskip.2cm\noindent{\bf Corollary \ref{EasyUniqueness}}
{\em For every fixed $k,m_1,\dots,m_k$ and family $\F$ of subsets of
$\{1,\dots,k+1\}$\break specifying a hierarchical collection of margins
to be disclosed, it can be decided in polynomial time whether any
specified entry $x_{i_1,\dots,i_{k+1}}$ is the same in all long
$m_1\times\cdots\times m_k\times n$ tables with the
disclosed margins, and hence at risk of exposure.}

\subsection{Terminology and Complexity}
\label{tac}

We use $\R$ for the reals, $\R_+$ for the nonnegative reals,
$\Z$ for the integers, and $\N$ for the nonnegative integers.
The sign of a real number $r$ is denoted by $\sign(r)\in\{0,-1,1\}$
and its absolute value is denoted by $|r|$. The $i$-th standard unit vector
in $\R^n$ is denoted by ${\bf 1}_i$. The {\em support} of $x\in\R^n$ is
the index set $\supp(x):=\{i:x_i\neq 0\}$ of nonzero entries of $x$.
The {\em indicator} of a subset $I\subseteq \{0,1\}^n$ is the vector
${\bf 1}_I:=\sum_{i\in I}{\bf 1}_i$ so that $\supp({\bf 1}_I)=I$.
When several vectors are indexed by subscripts, $w_1,\dots,w_d\in\R^n$,
their entries are indicated by pairs of subscripts,
$w_i=(w_{i,1},\dots, w_{i,n})$. When vectors are indexed
by superscripts, $x^1,\dots,x^k\in\R^n$, their entries are
indicated by subscripts, $x^i=(x^i_1,\dots, x^i_n)$.
The integer lattice $\Z^n$ is naturally embedded in $\R^n$.
The space $\R^n$ is endowed with the standard inner product which,
for $w,x\in\R^n$, is given by $wx:=\sum_{i=1}^n w_ix_i$.
Vectors $w$ in $\R^n$ will also be regarded as linear functionals
on $\R^n$ via the inner product $wx$. Thus, we refer to elements
of $\R^n$ as points, vectors, or linear functionals, as will be
appropriate from the context. The {\em convex hull} of a set $S\subseteq\R^n$
is denoted by $\conv(S)$ and the set of {\em vertices} of a polyhedron
$P\subseteq\R^n$ is denoted by $\vert(P)$. In linear discrete optimization
over $S\subseteq\Z^n$, the {\em facets} of $\conv(S)$ play an important
role, see Chv\'atal \cite{Chv} and the references therein for earlier work,
and Gr\"otschel, Lov\'asz and Schrijver \cite{GLS,Lov} for the later
culmination in the equivalence of separation and linear optimization
via the ellipsoid method of Yudin and Nemirovskii \cite{YN}.
As will turn out in \S \ref{rctldo}, in convex discrete
optimization over $S$, the {\em edges} of $\conv(S)$ play an
important role (most significantly in a way which is
{\em not} related to the Hirsch conjecture discussed
in \cite{KK}). We therefore use extensively convex polytopes,
for which we follow the terminology of \cite{Gru,Zie}.

We often assume that the feasible set $S\subseteq\Z^n$ is finite.
We then define its {\em radius} to be its $l_\infty$ radius
$\rho(S):=\max\{\|x\|_{\infty}:x\in S\}$ where, as usual,
$\|x\|_{\infty}:=\max_{i=1}^n |x_i|$. In other words, $\rho(S)$ is the
smallest $\rho\in\N$ such that $S$ is contained in the cube $[-\rho,\rho]^n$.

Our algorithms are applied to rational data only, and the time complexity
is as in the standard Turing machine model, see e.g. \cite{AHU,GJ,Sch}.
The input typically consists of rational (usually integer) numbers,
vectors, matrices, and finite sets of such objects.
The {\em binary length} of an integer number $z\in\Z$ is defined
to be the number of bits in its binary representation,
$\l z \r:=1+\lceil \log_2(|z|+1)\rceil$ (with the extra bit for the sign).
The length of a rational number presented as a fraction $r={p\over q}$
with $p,q\in\Z$ is $\l r \r:=\l p\r + \l q\r$. The length of an
$m\times n$ matrix $A$ (and in particular of a vector) is the sum
$\l A \r:=\sum_{i,j}\l a_{i,j}\r$ of the lengths of its entries.
Note that the length of $A$ is no smaller than the number of entries,
$\l A \r\geq mn$. Therefore, when $A$ is, say, part of an input to an
algorithm, with $m,n$ variable, the length $\l A \r$ already incorporates
$mn$, and so we will typically not account additionally for $m,n$ directly.
But sometimes, especially in results related to\break $n$-fold integer
programming, we will also emphasize $n$ as part of the input length.
Similarly, the length of a finite set $E$ of numbers, vectors or matrices
is the sum of lengths of its elements and hence, since $\l E\r\geq |E|$,
automatically accounts for its cardinality.

Some input numbers affect the running time of some algorithms through their
unary presentation, resulting in so-called ``pseudo polynomial" running time.
The {\em unary length} of an integer number $z\in\Z$ is the number $|z|+1$ of
bits in its unary representation (again, an extra bit for the sign). The unary
length of a rational number, vector, matrix, or finite set of such objects
are defined again as the sums of lengths of their numerical constituents,
and is again no smaller than the number of such numerical constituents.

When studying convex and linear integer programming in \S \ref{lnfip} and
\S \ref{cip} we sometimes have lower and upper bound vectors $l,u$ with
entries in $\Z_{\infty}:=\Z\uplus\{\pm\infty\}$. Both binary and unary
lengths of a $\pm\infty$ entry are constant, say $3$ by encoding
$\pm\infty:=\pm ``00"$.

To make the input encoding precise, we introduce the following notation.
In every algorithmic statement we describe explicitly the input encoding,
by listing in square brackets all input objects affecting the running time.
Unary encoded objects are listed directly whereas binary encoded objects
are listed in terms of their length. For example, as is often the case,
if the input of an algorithm consists of binary encoded vectors
(linear functionals) $w_1,\dots,w_d\in\Z^n$ and unary encoded integer
$\rho\in\N$ (bounding the radius $\rho(S)$ of the feasible set) then we will
indicate that the input is {\em encoded as} $[\rho,\l w_1,\dots,w_d\r]$.

Some of our algorithms are strongly polynomial time in the sense of \cite{Tar}.
For this, part of the input is regarded as ``special".
An algorithm is then {\em strongly polynomial time} if it is polynomial time
in the usual Turing sense with respect to all input, and in addition, the number
of arithmetic operations (additions, subtractions, multiplications, divisions,
and comparisons) it performs is polynomial in the special part of the input.
To make this precise, we extend our input encoding notation above by splitting
the square bracketed expression indicating the input encoding into a ``left"
side and a ``right" side, separated by semicolon, where the entire input is
described on the right and the special part of the input on the left.
For example, Theorem \ref{ConvexToLinear}, asserting that the algorithm
underlying it is strongly polynomial with data {\em encoded as}
$[n,|E|;\l\rho(S),w_1,\dots,w_d,E\r]$, where $\rho(S)\in\N$,
$w_1,\dots,w_d\in\Z^n$ and $E\subset\Z^n$, means that the running time
is polynomial in the binary length of $\rho(S)$, $w_1,\dots,w_d$,
and $E$, and the number of arithmetic operations is polynomial in $n$ and the
cardinality $|E|$, which constitute the special part of the input.

Often, as in \cite{GLS}, part of the input is presented by oracles. Then the
running time and the number of arithmetic operations count also the number of
oracle queries. An oracle algorithm is {\em polynomial time} if its running
time, including the number of oracle queries, and the manipulations of numbers,
some of which are answers to oracle queries, is polynomial in the length of
the input encoding. An oracle algorithm is {\em strongly polynomial time}
(with specified input encoding as above), if it is polynomial time in the
entire input (on the ``right"), and in addition, the number of arithmetic
operations it performs (including oracle queries) is
polynomial in the special part of the input (on the ``left").

\newpage
\section{Reducing Convex to Linear Discrete Optimization}
\label{rctldo}

In this section we show that when suitable auxiliary geometric information
about the convex hull $\conv(S)$ of a finite set $S\subseteq\Z^n$ is available,
the convex discrete optimization problem over $S$ can be reduced to the
solution of strongly polynomially many linear discrete optimization counterparts
over $S$. This result will be incorporated into the polynomial time algorithms
developed in \S \ref{ccoam} and \S \ref{cip} for convex combinatorial
optimization and convex integer programming respectively.
In \S \ref{edaz} we provide some preliminaries on edge-directions and
zonotopes. In \S \ref{sproctldo} we prove the reduction which is the main
result of this section. In \S \ref{pprwedana} we prove a pseudo polynomial
reduction for any finite set.

\subsection{Edge-Directions and Zonotopes}
\label{edaz}

We begin with some terminology and facts that
play an important role in the sequel.
A {\em direction} of an edge ($1$-dimensional face) $e=[u,v]$ of a polytope
$P$ is any nonzero scalar multiple of $u-v$. A set of vectors $E$
{\em covers all edge-directions of $P$} if it contains a direction of
each edge of $P$.
The {\em normal cone} of a polytope $P\subset\R^n$ at its
face $F$ is the (relatively open) cone $C_P^F$ of those linear functionals
$h\in\R^n$ which are maximized over $P$ precisely at points of $F$.
A polytope $Z$ is a {\em refinement} of a polytope $P$ if the normal cone
of every vertex of $Z$ is contained in the normal cone of some vertex of $P$.
If $Z$ refines $P$ then, moreover, the closure of each
normal cone of $P$ is the union of closures of normal cones of $Z$.
The {\em zonotope} generated by a set of vectors
$E=\{e_1,\dots,e_m\}$ in $\R^d$ is the following polytope,
which is the projection by $E$ of the cube $[-1,1]^m$ into $\R^d$,
$$Z\quad:=\quad\zone(E)\quad:=\quad
\conv\left\{\sum_{i=1}^m \lambda_i e_i\,:\,\lambda_i=\pm 1\right\}
\quad\subset\quad\R^d\quad .$$

The following fact goes back to Minkowski, see \cite{Gru}.
\bl{Refinement}
Let $P$ be a polytope and let $E$ be a finite set that covers
all edge-directions of $P$. Then the zonotope $Z:=\zone(E)$
generated by $E$ is a refinement of $P$.
\el
\bpr
Consider any vertex $u$ of $Z$. Then $u=\sum_{e\in E} \lambda_e e$ for
suitable $\lambda_e=\pm1$. Thus, the normal cone $C^u_Z$ consists of
those $h$ satisfying $h\lambda_e e > 0$ for all $e$. Pick any
${\hat h}\in C^u_Z$ and let $v$ be a vertex of $P$ at which $\hat h$ is
maximized over $P$. Consider any edge $[v,w]$ of $P$.
Then $v-w=\alpha_e e$ for some scalar $\alpha_e\neq 0$ and some $e\in E$,
and $0\leq {\hat h}(v-w)={\hat h}\alpha_e e$, implying $\alpha_e\lambda_e>0$.
It follows that every $h\in C^u_Z$ satisfies $h(v-w)>0$ for every edge of $P$
containing $v$. Therefore $h$ is maximized over $P$ uniquely at $v$
and hence is in the cone $C^v_P$ of $P$ at $v$.
This shows $C^u_Z\subseteq C^v_P$. Since $u$ was arbitrary, it follows
that the normal cone of every vertex of $Z$ is contained in the normal cone
of some vertex of $P$.
\epr
The next lemma provides bounds on the number of vertices of any zonotope
and on the algorithmic complexity of constructing its vertices,
each vertex along with a linear functional maximized over the
zonotope uniquely at that vertex. The bound on the number of vertices
has been rediscovered many times over the years. An early reference
is \cite{Har}, stated in the dual form of $2$-partitions.
A more general treatment is \cite{Zas}. Recent extensions
to $p$-partitions for any $p$ are in \cite{AO,HOR}, and to Minkowski sums
of arbitrary polytopes are in \cite{GS}. Interestingly, already in \cite{Har},
back in 1967, the question was raised about the algorithmic
complexity of the problem; this is now settled in \cite{EOS,ESS}
(the latter reference correcting the former). We state the precise bounds
on the number of vertices and arithmetic complexity, but will need later
only that for any fixed $d$ the bounds are polynomial in the
number of generators. Therefore, below we only outline a proof that the
bounds are polynomial. Complete details are in the above references.
\bl{Zonotope}
The number of vertices of any zonotope $Z:=\zone(E)$ generated by a set
$E$ of $m$ vectors in $\R^d$ is at most $2\sum_{k=0}^{d-1}{{m-1}\choose k}$.
For every fixed $d$, there is a strongly polynomial time algorithm that,
given $E\subset\Z^d$, encoded as $[m:=|E|;\l E \r]$, outputs every vertex $v$
of $Z:=\zone(E)$ along with a linear functional $h_v\in\Z^d$ maximized over $Z$
uniquely at $v$, using $O(m^{d-1})$ arithmetics operations for $d\geq 3$
and $O(m^d)$ for $d\leq 2$.
\el
\bpr
We only outline a proof that, for every fixed $d$, the polynomial bounds
$O(m^{d-1})$ on the number of vertices and $O(m^{d})$ on the arithmetic
complexity hold. We assume that $E$ linearly spans $\R^d$ (else the dimension
can be reduced) and is generic, that is, no $d$ points of $E$ lie on a
linear hyperplane (one containing the origin). In particular, $0\notin E$.
The same bound for arbitrary $E$ then follows using a perturbation argument
(cf. \cite{HOR}).

Each oriented linear hyperplane $H=\{x\in\R^d:hx=0\}$ with $h\in\R^d$
nonzero induces a partition of $E$ by $E=H^-\biguplus H^0\biguplus H^+$,
with $H^-:=\{e\in E:he<0\}$, $E^0:=E\cap H$, and $H^+:=\{e\in E:he>0\}$.
The vertices of $Z=\zone(E)$ are in bijection with ordered $2$-partitions
of $E$ induced by such hyperplanes that avoid $E$.
Indeed, if $E=H^-\biguplus H^+$ then the linear functional $h_v:=h$
defining $H$ is maximized over $Z$ uniquely at the vertex
$v:=\sum\{e:e\in H^+\}-\sum\{e:e\in H^-\}$ of $Z$.

We now show how to enumerate all such $2$-partitions and hence vertices
of $Z$. Let $M$ be any of the $m\choose d-1$ subsets of $E$ of size $d-1$.
Since $E$ is generic, $M$ is linearly independent and spans a unique
linear hyperplane $\lin(M)$. Let ${\hat H}=\{x\in\R^d:{\hat h}x=0\}$ be
one of the two orientations of the hyperplane $\lin(M)$. Note that
${\hat H}^0=M$. Finally, let $L$ be any of the $2^{d-1}$ subsets of $M$.
Since $M$ is linearly independent, there is a $g\in\R^d$ which linearly
separates $L$ from $M\setminus L$, namely, satisfies $gx<0$ for all
$x\in L$ and $gx>0$ for all $x\in M\setminus L$. Furthermore, there is a
sufficiently small $\epsilon>0$ such that the oriented hyperplane
$H:=\{x\in\R^d:hx=0\}$ defined by $h:={\hat h}+\epsilon g$ avoids $E$
and the $2$-partition induced by $H$ satisfies $H^-={\hat H}^-\biguplus L$
and $H^+={\hat H}^+\biguplus(M\setminus L)$. The corresponding
vertex of $Z$ is $v:=\sum\{e:e\in H^+\}-\sum\{e:e\in H^-\}$ and
the corresponding linear functional which is maximized
over $Z$ uniquely at $v$ is $h_v:=h={\hat h}+\epsilon g$.

We claim that any ordered $2$-partition arises that way from some $M$,
some orientation ${\hat H}$ of $\lin(M)$, and some $L$. Indeed, consider
any oriented linear hyperplane $\tilde H$ avoiding $E$. It can
be perturbed to a suitable oriented $\hat H$ that touches precisely
$d-1$ points of $E$. Put $M:={\hat H}^0$ so that $\hat H$ coincides
with one of the two orientations of the hyperplane $\lin(M)$ spanned by $M$,
and put $L:={\tilde H}^-\cap M$. Let $H$ be an oriented hyperplane obtained
from $M$, $\hat H$ and $L$ by the above procedure. Then the ordered
$2$-partition $E=H^-\biguplus H^+$ induced by $H$ coincides with the ordered
$2$-partition $E={\tilde H}^-\biguplus {\tilde H}^+$ induced by $\tilde H$.

Since there are $m\choose d-1$ many $(d-1)$-subsets $M\subseteq E$, two
orientations ${\hat H}$ of $\lin(M)$, and $2^{d-1}$ subsets $L\subseteq M$,
and $d$ is fixed, the total number of $2$-partitions and hence also the total
number of vertices of $Z$ obey the upper bound $2^d{m\choose d-1}=O(m^{d-1})$.
Furthermore, for each choice of $M$, ${\hat H}$ and $L$,
the linear functional $\hat h$ defining ${\hat H}$, as well as $g$,
$\epsilon$, $h_v=h={\hat h}+\epsilon g$, and the vertex
$v=\sum\{e:e\in H^+\}-\sum\{e:e\in H^-\}$ of $Z$ at which $h_v$ is uniquely
maximized over $Z$, can all be computed using $O(m)$ arithmetic operations.
This shows the claimed bound $O(m^d)$ on the arithmetic complexity.
\epr

We conclude with a simple fact about
edge-directions of projections of polytopes.
\bl{Projection}
If $E$ covers all edge-directions of a polytope $P$, and $Q:=\omega(P)$
is the image of $P$ under a linear map $\omega:\R^n\longrightarrow \R^d$,
then $\omega(E)$ covers all edge-directions of $Q$.
\el
\bpr
Let $f$ be a direction of an edge $[x,y]$ of $Q$. Consider the face
$F:=\omega^{-1}([x,y])$ of $P$. Let $V$ be the set of vertices of $F$
and let $U=\{u\in V\,:\,\omega(u)=x\,\}$. Then for some $u\in U$ and
$v\in V\setminus U$, there must be an edge $[u,v]$ of $F$, and hence
of $P$. Then $\omega(v)\in(x,y]$ hence $\omega(v)=x+\alpha f$ for
some $\alpha\neq 0$. Therefore, with $e:={1\over\alpha}(v-u)$,
a direction of the edge $[u,v]$ of $P$, we find that
$f={1\over\alpha}(\omega(v)-\omega(u))=\omega(e)\in\omega(E)$.
\epr

\subsection{Strongly Polynomial Reduction
of Convex to Linear Discrete Optimization}
\label{sproctldo}

A {\em linear discrete optimization oracle} for a set $S\subseteq\Z^n$
is one that, queried on $w\in\Z^n$, either returns an optimal solution to
the linear discrete optimization problem over $S$, that is, an $x^*\in S$
satisfying $wx^*=\max\{wx:x\in S\}$, or asserts that none exists, that is,
either the problem is infeasible or the objective function is unbounded.
We now show that a set $E$ covering all edge-directions of the polytope
$\conv(S)$ underlying a convex discrete optimization problem over a finite
set $S\subset\Z^n$ allows to solve it by solving polynomially many linear
discrete optimization counterparts over $S$. The following theorem extends
and unifies the corresponding reductions in \cite{OR} and \cite{DHORW} for
convex combinatorial optimization and convex integer programming respectively.
Recall from \S \ref{tac} that the {\em radius} of a finite set
$S\subset\Z^n$ is defined to be $\rho(S):=\max\{|x_i|:x\in S,\ i=1,\dots,n\}$.

\bt{ConvexToLinear}
For every fixed $d$ there is a strongly polynomial time algorithm that,
given finite set $S\subset\Z^n$ presented by a linear discrete optimization
oracle, integer vectors $w_1,\dots,w_d\in\Z^n$, set $E\subset\Z^n$ covering
all edge-directions of $\conv(S)$, and convex functional
$c:\R^d\longrightarrow\R$ presented by a comparison oracle, encoded as
$[n,|E|;\l \rho(S),w_1,\dots,w_d,E\r]$,
solves the convex discrete optimization problem
$$\max\, \{c(w_1 x,\dots,w_d x):\ x\in S\}\ .$$
\et
\bpr
First, query the linear discrete optimization oracle presenting $S$ on
the trivial linear functional $w=0$. If the oracle asserts that there
is no optimal solution then $S$ is empty so terminate the algorithm
asserting that no optimal solution exists to the convex discrete
optimization problem either. So assume the problem is feasible.
Let $P:=\conv(S)\subset\R^n$ and $Q:=\{(w_1x,\dots,w_dx):x\in P\}\subset\R^d$.
Then $Q$ is a projection of $P$, and hence by Lemma \ref{Projection}
the projection $D:=\{(w_1e,\dots,w_de):e\in E\}$ of the set $E$
is a set covering all edge-directions of $Q$. Let $Z:=\zone(D)\subset\R^d$
be the zonotope generated by $D$. Since $d$ is fixed, by Lemma \ref{Zonotope}
we can produce in strongly polynomial time all vertices of $Z$, every
vertex $v$ along with a linear functional $h_v\in\Z^d$ maximized
over $Z$ uniquely at $v$. For each of these polynomially many $h_v$,
repeat the following procedure. Define a vector $g_v\in\Z^n$
by $g_{v,j}:=\sum_{i=1}^d w_{i,j}h_{v,i}$ for $j=1,\dots,n$.
Now query the linear discrete optimization oracle presenting $S$ on the
linear functional $w:=g_v\in\Z^n$. Let $x_v\in S$ be the optimal solution
obtained from the oracle, and let $z_v:=(w_1x_v,\dots,w_dx_v)\in Q$
be its projection. Since $P=\conv(S)$, we have that $x_v$ is also
a maximizer of $g_v$ over $P$. Since for every $x\in P$ and its projection
$z:=(w_1x,\dots,w_dx)\in Q$ we have $h_vz=g_vx$, we conclude that
$z_v$ is a maximizer of $h_v$ over $Q$. Now we claim that each vertex
$u$ of $Q$ equals some $z_v$. Indeed, since $Z$ is a refinement of $Q$
by Lemma \ref{Refinement}, it follows that there is some vertex $v$
of $Z$ such that $h_v$ is maximized over $Q$ uniquely at $u$, and
therefore $u=z_v$. Since $c(w_1x,\dots,w_dx)$ is convex on
$\R^n$ and $c$ is convex on $\R^d$, we find that
\begin{eqnarray*}
\hskip-0.8cm\max_{x\in S}\, c(w_1 x,\dots,w_d x)
& \, =\,  & \max_{x\in P}\, c(w_1 x,\dots,w_d x)
\quad\quad =\, \max_{z\in Q}\, c(z) \\
& = & \max\{c(u):u\ \mbox{vertex of }Q\}
\, =\, \max\{c(z_v): v\ \mbox{vertex of } Z\}\,.
\end{eqnarray*}
Using the comparison oracle of $c$, find a vertex $v$ of $Z$
attaining maximum value $c(z_v)$, and output $x_v\in S$,
an optimal solution to the convex discrete optimization problem.
\epr

\subsection{Pseudo Polynomial Reduction when Edge-Directions are not Available}
\label{pprwedana}

Theorem \ref{ConvexToLinear} reduces convex discrete optimization
to polynomially many linear discrete optimization counterparts when a
set covering all edge-directions of the underlying polytope is available.
However, often such a set is not available
(see e.g. \cite{BO} for the important case of bipartite matching).
We now show how to reduce convex discrete optimization to many linear
discrete optimization counterparts when a set covering all edge-directions
is not offhand available. In the absence of such a set, the problem is much
harder, and the algorithm below is polynomially bounded only in the unary
length of the radius $\rho(S)$ and of the linear functionals $w_1,\dots,w_d$,
rather than in their binary length $\l\rho(S),w_1,\dots, w_d\r$
as in the algorithm of Theorem \ref{ConvexToLinear}.
Moreover, an upper bound $\rho\geq\rho(S)$ on the radius of $S$ is
required to be given explicitly in advance as part of the input.

\bt{NoEdgeDirections}
For every fixed $d$ there is a polynomial time algorithm that, given finite set
$S\subseteq\Z^n$ presented by a linear discrete optimization oracle, integer
$\rho\geq \rho(S)$, vectors $w_1,\dots,w_d\in\Z^n$, and convex functional
$c:\R^d\longrightarrow\R$ presented by a comparison oracle, encoded as
$[\rho, w_1,\dots,w_d]$, solves the convex discrete optimization problem
$$\max\, \{c(w_1 x,\dots,w_d x):\ x\in S\}\ .$$
\et

\bpr
Let $P:=\conv(S)\subset\R^n$, let $T:=\{(w_1x,\dots,w_dx):x\in S\}$
be the projection of $S$ by $w_1,\dots,w_d$, and let
$Q:=\conv(T)\subset\R^d$ be the corresponding projection of $P$.
Let $r:=n\rho\max_{i=1}^d\!\|w_i\|_{\infty}$ and let
$G:=\{-r,\dots,-1,0,1,\dots,r\}^d$. Then $T\subseteq G$ and the number
$(2r+1)^d$ of points of $G$ is polynomially bounded in the input as encoded.

Let $D:=\{u-v:u,v\in G, u\neq v\}$ be the set of differences of pairs
of distinct point of $G$. It covers all edge-directions of $Q$ since
$\vert(Q)\subseteq T\subseteq G$. Moreover, the number of points of $D$
is less than $(2r+1)^{2d}$ and hence polynomial in the input.
Now invoke the algorithm of Theorem \ref{ConvexToLinear}: while the
algorithm requires a set $E$ covering all edge-directions of $P$, it needs $E$
only to compute a set $D$ covering all edge-directions of the projection $Q$
(see proof of Theorem \ref{ConvexToLinear}), which here is computed directly.
\epr

\newpage
\section{Convex Combinatorial Optimization and More}
\label{ccoam}

In this section we discuss convex combinatorial optimization.
The main result is that convex combinatorial optimization over
a set $S\subseteq\{0,1\}^n$ presented by a membership oracle
can be solved in strongly polynomial time provided it is endowed
with a set covering all edge-directions of $\conv(S)$.
In particular, the standard linear combinatorial optimization problem over
$S$ can be solved in strongly polynomial time as well. In \S \ref{fmtlo}
we provide some preparatory statements involving various oracle presentation
of the feasible set $S$. In \S \ref{laccoispt} we combine these preparatory
statements with Theorem \ref{ConvexToLinear} and prove the main result of
this section. An extension to arbitrary finite sets $S\subset\Z^n$
endowed with edge-directions is established in \S \ref{lacdooasippt}.
We conclude with some applications in \S \ref{sacco}.

As noted in the introduction, when $S$ is contained in $\{0,1\}^n$ we refer
to discrete optimization over $S$ also as {\em combinatorial optimization}
over $S$, to emphasize that $S$ typically represents a family
$\F\subseteq 2^N$ of subsets of a ground set $N:=\{1,\dots,n\}$
possessing some combinatorial property of interest
(for instance, the family of bases of a matroid over $N$,
see \S \ref{mamnst}). The convex combinatorial optimization problem
then also has the following interpretation (taken in \cite{Onn1,OR}).
We are given a weighting\break $\omega:N\longrightarrow\Z^d$
of elements of the ground set by $d$-dimensional integer vectors.
We interpret the weight vector $\omega(j)\in\Z^d$ of element $j$ as
representing its value under $d$ criteria (e.g., if $N$ is the set
of edges in a network then such criteria may include profit,
reliability, flow velocity, etc.). The weight of a subset
$F\subseteq N$ is the sum $\omega(F):=\sum_{j\in F}\omega(j)$
of weights of its elements, representing the total value of $F$ under
the $d$ criteria. Now, given a convex functional $c:\R^d\longrightarrow\R$,
the objective function value of $F\subseteq N$ is the ``convex balancing"
$c(\omega(F))$ of the values of the weight vector of $F$.
The convex combinatorial optimization problem is to find a family member
$F\in\F$ maximizing $c(\omega(F))$. The usual linear combinatorial
optimization problem over $\F$ is the special case of $d=1$ and
$c$ the identity on $\R$. To cast a problem of that form in our usual
setup just let $S:=\{{\bf 1}_F:F\in\F\}\subseteq\{0,1\}^n$ be the set
of indicators of members of $\F$ and define weight vectors
$w_1,\dots,w_d\in\Z^n$ by $w_{i,j}:=\omega(j)_i$ for
$i=1,\dots,d$ and $j=1,\dots,n$.

\subsection{From Membership to Linear Optimization}
\label{fmtlo}

A {\em membership oracle} for a set $S\subseteq\Z^n$ is
one that, queried on $x\in\Z^n$, asserts whether or not $x\in S$.
An {\em augmentation oracle} for $S$ is one that, queried on $x\in S$ and
$w\in\Z^n$, either returns an ${\hat x}\in S$ with $w{\hat x}>wx$, i.e. a
better point of $S$, or asserts that none exists, i.e. $x$ is
optimal for the linear discrete optimization problem over $S$.

A membership oracle presentation of $S$ is very broad and available in all
reasonable applications, but reveals little information on $S$, making it
hard to use. However, as we now show, the edge-directions of $\conv(S)$
allow to convert membership to augmentation.

\bl{Membership}
There is a strongly polynomial time algorithm that, given set
$S\subseteq\{0,1\}^n$ presented by a membership oracle, $x\in S$,
$w\in\Z^n$, and set $E\subset\Z^n$ covering all edge-directions of
the polytope $\conv(S)$, encoded as $[n,|E|;\l x,w,E \r]$,
either returns a better point ${\hat x}\in S$, that is,
one satisfying $w{\hat x}>wx$, or asserts that none exists.
\el

\bpr
Each edge of $P:=\conv(S)$ is the difference of two $\{0,1\}$-vectors.
Therefore, each edge direction of $P$ is, up to scaling,
a $\{-1,0,1\}$-vector. Thus, scaling $e:={1\over{\|e\|_{\infty}}}e$
and $e:=-e$ if necessary, we may and will assume that $e\in\{-1,0,1\}^n$
and $we\geq 0$ for all $e\in E$. Now, using the membership oracle,
check if there is an $e\in E$ such that $x+e\in S$ and $we>0$. If there is
such an $e$ then output ${\hat x}:=x+e$ which is a better point, whereas if
there is no such $e$ then terminate asserting that no better point exists.

Clearly, if the algorithm outputs an $\hat x$ then it is indeed a better
point. Conversely, suppose $x$ is not a maximizer of $w$ over $S$.
Since $S\subseteq\{0,1\}^n$, the point $x$ is a vertex of $P$. Since $x$
is not a maximizer of $w$, there is an edge $[x,{\hat x}]$ of $P$ with
$\hat x$ a vertex satisfying $w{\hat x}>wx$. But then $e:={\hat x}-x$ is the
one $\{-1,0,1\}$ edge-direction of $[x,{\hat x}]$ with $we\geq 0$ and hence
$e\in E$. Thus, the algorithm will find and output ${\hat x}=x+e$ as it should.
\epr

An augmentation oracle presentation of a finite $S$ allows to solve the linear
discrete optimization problem $\max\{wx:x\in S\}$ over $S$ by starting from
any feasible $x\in S$ and repeatedly augmenting it until an optimal solution
$x^*\in S$ is reached. The next lemma bounds the running time needed to
reach optimality using this procedure. While the running time is polynomial
in the binary length of the linear functional $w$ and the initial point $x$,
it is more sensitive to the radius $\rho(S)$ of the feasible set $S$,
and is polynomial only in its unary length. The lemma is an adaptation
of a result of \cite{GL,SWZ} (stated therein for $\{0,1\}$-sets),
which makes use of bit-scaling ideas going back to \cite{EK}.

\bl{Augmentation}
There is a polynomial time algorithm that, given finite set
$S\subset\Z^n$ presented by an augmentation oracle, $x\in S$, and
$w\in\Z^n$, encoded as $[\rho(S),\l x,w \r]$, provides an optimal solution
$x^*\in S$ to the linear discrete optimization problem $\max\{wz:z\in S\}$.
\el
\bpr
Let $k:=\max_{j=1}^n\lceil\log_2 (|w_j|+1)\rceil$ and note that $k\leq\l w\r$.
For $i=0,\dots,k$ define a linear functional $u_i=(u_{i,1},\dots,u_{i,n})\in\Z^n$
by $u_{i,j}:=\sign(w_j)\lfloor 2^{i-k}|w_j|\rfloor$ for $j=1,\dots,n$. Then
$u_0=0$, $u_k=w$, and $u_i-2u_{i-1}\in\{-1,0,1\}^n$ for all $i=1,\dots,k$.

We now describe how to construct a sequence of points $y_0,y_1,\dots,y_k\in S$
such that $y_i$ is an optimal solution to $\max\{u_iy:y\in S\}$ for all $i$.
First note that all points of $S$ are optimal for $u_0=0$  and hence we can
take $y_0:=x$ to be the point of $S$ given as part of the input.
We now explain how to determine $y_i$ from $y_{i-1}$ for $i=1,\dots,k$.
Suppose $y_{i-1}$ has been determined. Set ${\tilde y}:=y_{i-1}$.
Query the augmentation oracle on ${\tilde y}\in S$ and $u_i$; if the oracle
returns a better point $\hat y$ then set ${\tilde y}:={\hat y}$ and repeat,
whereas if it asserts that there is no better point then the optimal solution
for $u_i$ is read off to be $y_i:={\tilde y}$. We now bound the number
of calls to the oracle. Each time the oracle is queried on ${\tilde y}$
and $u_i$ and returns a better point ${\hat y}$, the improvement is by at
least one, i.e. $u_i({\hat y}-{\tilde y})\geq 1$; this is so
because $u_i$, ${\tilde y}$ and ${\hat y}$ are integer. Thus, the number
of necessary augmentations from $y_{i-1}$ to $y_i$ is at most
the total improvement, which we claim satisfies
$$u_i(y_i-y_{i-1})
=(u_i-2u_{i-1})(y_i-y_{i-1})+2u_{i-1}(y_i-y_{i-1})\leq 2n\rho+0=2n\rho\ ,$$
where $\rho:=\rho(S)$. Indeed, $u_i-2u_{i-1}\in\{-1,0,1\}^n$ and
$y_i,y_{i-1}\in S\subset [-\rho,\rho]^n$ imply
$(u_i-2u_{i-1})(y_i-y_{i-1})\leq 2n\rho$; and $y_{i-1}$
optimal for $u_{i-1}$ gives $u_{i-1}(y_i-y_{i-1})\leq 0$.

Thus, after a total number of at most $2n\rho k$ calls to the oracle we obtain
$y_k$ which is optimal for $u_k$. Since $w=u_k$ we can output $x^*:=y_k$ as the
desired optimal solution to the linear discrete optimization problem. Clearly
the number $2n\rho k$ of calls to the oracle, as well as the number of
arithmetic operations and binary length of numbers occurring during the
algorithm, are polynomial in $\rho(S),\l x,w \r$. This completes the proof.
\epr

We conclude this preparatory subsection by recording the following
result of \cite{FT} which incorporates the heavy simultaneous
Diophantine approximation of \cite{LLL}.

\bp{Diophantine}
There is a strongly polynomial time algorithm that, given $w\in\Z^n$,
encoded as $[n;\l w \r]$, produces ${\hat w}\in\Z^n$, whose binary length
$\l{\hat w}\r$ is polynomially bounded in $n$ and independent of $w$,
and with $\sign({\hat w}z)=\sign(wz)$ for every $z\in\{-1,0,1\}^n$.
\ep

\subsection{Linear and Convex Combinatorial
Optimization in Strongly Polynomial Time}
\label{laccoispt}

Combining the preparatory statements of \S \ref{fmtlo} with
Theorem \ref{ConvexToLinear}, we can now solve the convex combinatorial
optimization over a set $S\subseteq\{0,1\}^n$ presented by a
membership oracle and endowed with a set covering all edge-directions
of $\conv(S)$ in strongly polynomial time.
We start with the special case of linear combinatorial optimization.

\bt{LCO}
There is a strongly polynomial time algorithm that, given set
$S\subseteq\{0,1\}^n$ presented by a membership oracle, $x\in S$, $w\in\Z^n$,
and set $E\subset\Z^n$ covering all edge-directions of the polytope $\conv(S)$,
encoded as $[n,|E|;\l x,w,E\r]$, provides an optimal solution $x^*\in S$
to the linear combinatorial optimization problem $\max\{wz:z\in S\}$.
\et

\bpr
First, an augmentation oracle for $S$ can be simulated using
the membership oracle, in strongly polynomial time,
by applying the algorithm of Lemma \ref{Membership}

Next, using the simulated augmentation oracle for $S$, we can now do linear
optimization over $S$ in strongly polynomial time as follows.
First, apply to $w$ the algorithm of
Proposition \ref{Diophantine} and obtain $\hat w\in\Z^n$ whose binary
length $\l{\hat w}\r$ is polynomially bounded in $n$, which satisfies
$\sign({\hat w}z)=\sign(wz)$ for every $z\in\{-1,0,1\}^n$.
Since $S\subseteq\{0,1\}^n$, it is finite and has radius $\rho(S)=1$.
Now apply the algorithm of Lemma \ref{Augmentation}
to $S$, $x$ and $\hat w$, and obtain a maximizer
$x^*$ of $\hat w$ over $S$. For every $y\in\{0,1\}^n$ we then have
$x^*-y\in\{-1,0,1\}^n$ and hence $\sign(w(x^*-y))=\sign({\hat w}(x^*-y))$.
So $x^*$ is also a maximizer of $w$ over $S$ and hence an
optimal solution to the given linear combinatorial optimization problem.
Now, $\rho(S)=1$, $\l{\hat w}\r$ is polynomial in $n$, and $x\in\{0,1\}^n$
and hence $\l x \r$ is linear in $n$. Thus, the entire length of the input
$[\rho(S),\l x,{\hat w}\r]$ to the polynomial-time algorithm of
Lemma \ref{Augmentation} is polynomial in $n$, and so its running time
is in fact strongly polynomial on that input.
\epr

Combining Theorems \ref{ConvexToLinear} and \ref{LCO}
we recover at once the following result of \cite{OR}.

\bt{CCO}
For every fixed $d$ there is a strongly polynomial time
algorithm that, given set $S\subseteq\{0,1\}^n$ presented by
a membership oracle, $x\in S$, vectors $w_1,\dots,w_d\in\Z^n$,
set $E\subset\Z^n$ covering all edge-directions of the polytope
$\conv(S)$, and convex functional $c:\R^d\longrightarrow\R$ presented
by a comparison oracle, encoded as $[n,|E|;\l x,w_1,\dots,w_d,E\r]$,
provides an optimal solution $x^*\in S$ to
the convex combinatorial optimization problem
$$\max\, \{c(w_1 z,\dots,w_d z):\ z\in S\}\ .$$
\et
\bpr
Since $S$ is nonempty, a linear discrete optimization oracle for $S$ can be
simulated in strongly polynomial time by the algorithm of Theorem \ref{LCO}.
Using this simulated oracle, we can apply
the algorithm of Theorem \ref{ConvexToLinear} and solve the given
convex combinatorial optimization problem in strongly polynomial time.
\epr

\subsection{Linear and Convex Discrete Optimization over
any Set in Pseudo Polynomial Time}
\label{lacdooasippt}

In \S \ref{laccoispt} above we developed strongly polynomial time
algorithms for linear and convex discrete optimization over $\{0,1\}$-sets.
We now provide extensions of these algorithms to arbitrary finite
sets $S\subset\Z^n$. As can be expected, the algorithms become slower.

We start by recording the following fundamental result of
Khachiyan \cite{Kha} asserting that linear programming is polynomial
time solvable via the ellipsoid method \cite{YN}. This result will be
used below as well as several more times later in the monograph.

\bp{LinearProgramming}
There is a polynomial time algorithm that, given $A\in\Z^{m\times n}$,
$b\in\Z^m$, and $w\in\Z^n$, encoded as $[\l A,b,w \r]$, either asserts that
$P:=\{x\in\R^n:Ax\leq b\}$ is empty, or asserts that the linear functional
$wx$ is unbounded over $P$, or provides a vertex $v\in\vert(P)$ which
is an optimal solution to the linear program $\max\{wx:x\in P\}$.
\ep

The following analog of Lemma \ref{Membership} shows how to covert
membership to augmentation in polynomial time, albeit, no longer in
strongly polynomial time. Here, both the given initial point $x$ and the
returned better point $\hat x$ if any, are {\em vertices} of $\conv(S)$.

\bl{MoreMembership}
There is a polynomial time algorithm that, given finite set $S\subset\Z^n$
presented by a membership oracle, vertex $x$ of the polytope $\conv(S)$,
$w\in\Z^n$, and set $E\subset\Z^n$ covering all edge-directions of $\conv(S)$,
encoded as $[\rho(S),\l x,w,E \r]$, either returns a better vertex ${\hat x}$ of
$\conv(S)$, that is, one satisfying $w{\hat x}>wx$, or asserts that none exists.
\el

\bpr
Dividing each vector $e\in E$ by the greatest common divisor of its entries
and setting $e:=-e$ if necessary, we can and will assume that each $e$ is
{\em primitive}, that is, its entries are relatively prime integers, and
$we\geq 0$. Using the membership oracle, construct the subset $F\subseteq E$
of those $e\in E$ for which $x+re\in S$ for some $r\in\{1,\dots,2\rho(S)\}$.
Let $G\subseteq F$ be the subset of those $f\in F$ for which $wf>0$.
If $G$ is empty then terminate asserting that there is no better vertex.
Otherwise, consider the convex cone $\cone(F)$ generated by $F$.
It is clear that $x$ is incident on an edge of $\conv(S)$ in
direction $f$ if and only if $f$ is an extreme ray of $\cone(F)$.
Moreover, since $G=\{f\in F:wf>0\}$ is nonempty, there must
be an extreme ray of $\cone(F)$ which lies in $G$.
Now $f\in F$ is an extreme ray of $\cone(F)$ if and only if
there do not exist nonnegative $\lambda_e$, $e\in F\setminus\{f\}$, such that
$f=\sum_{e\neq f}\lambda_e e$; this can be checked in polynomial time using
linear programming. Applying this procedure to each $f\in G$,
identify an extreme ray $g\in G$. Now, using the membership oracle,
determine the largest $r\in\{1,\dots,2\rho(S)\}$ for which $x+rg\in S$.
Output ${\hat x}:=x+rg$ which is a better vertex of $\conv(S)$.
\epr

We now prove the extensions of Theorems \ref{LCO} and \ref{CCO} to
arbitrary, not necessarily $\{0,1\}$-valued, finite sets. While the running
time remains polynomial in the binary length of the weights $w_1,\dots,w_d$
and the set of edge-directions $E$, it is more sensitive to
the radius $\rho(S)$ of the feasible set $S$, and is polynomial only
in its unary length. Here, the initial feasible point and the optimal
solution output by the algorithms are vertices of $\conv(S)$.
Again, we start with the special case of linear combinatorial optimization.

\bt{LinearArbitrarySets}
There is a polynomial time algorithm that, given finite $S\subset\Z^n$
presented by a membership oracle, vertex $x$ of the polytope $\conv(S)$, $w\in\Z^n$,
and set $E\subset\Z^n$ covering all edge-directions of $\conv(S)$, encoded
as $[\rho(S),\l x,w,E\r]$, provides an optimal solution $x^*\in S$
to the linear discrete optimization problem $\max\{wz:z\in S\}$.
\et
\bpr
Apply the algorithm of Lemma \ref{Augmentation} to the given data.
Consider any query $x'\in S$, $w'\in\Z^n$ made by that algorithm to an
augmentation oracle for $S$. To answer it, apply the algorithm of
Lemma \ref{MoreMembership} to $x'$ and $w'$. Since the first query made by
the algorithm of Lemma \ref{Augmentation} is on the given input
vertex $x':=x$, and any consequent query is on a point $x':=\hat x$
which was the reply of the augmentation oracle to the previous query
(see proof of Lemma \ref{Augmentation}), we see that the algorithm of
Lemma \ref{MoreMembership} will always be asked on a vertex of $S$ and reply
with another. Thus, the algorithm of Lemma \ref{MoreMembership}
can answer all augmentation queries and enables the polynomial
time solution of the given problem.
\epr

\bt{ConvexArbitrarySets}
For every fixed $d$ there is a polynomial time algorithm that, given finite
set $S\subseteq\Z^n$ presented by membership oracle, vertex $x$ of $\conv(S)$,
vectors $w_1,\dots,w_d\in\Z^n$, set $E\subset\Z^n$ covering all edge-directions
of the polytope $\conv(S)$, and convex functional $c:\R^d\longrightarrow\R$
presented by a comparison oracle, encoded as $[\rho(S),\l x,w_1,\dots,w_d,E\r]$,
provides an optimal solution $x^*\in S$ to
the convex combinatorial optimization problem
$$\max\, \{c(w_1 z,\dots,w_d z):\ z\in S\}\ .$$
\et

\bpr
Since $S$ is nonempty, a linear discrete optimization oracle for $S$
can be simulated in polynomial time by the algorithm of
Theorem \ref{LinearArbitrarySets}. Using this simulated oracle,
we can apply the algorithm of Theorem \ref{ConvexToLinear}
and solve the given problem in polynomial time.
\epr

\subsection{Some Applications}
\label{sacco}

\subsubsection{Positive Semidefinite Quadratic Binary Programming}
\label{psqbp}

The quadratic binary programming problem is the following: given an $n\times n$
matrix $M$, find a vector $x\in\{0,1\}^n$ maximizing the quadratic form
$x^TMx$ induced by $M$. We consider here the instance where $M$ is positive
semidefinite, in which case it can be assumed to be presented as $M=W^TW$
with $W$ a given $d\times n$ matrix. Already this restricted version is
very broad: if the rank $d$ of $W$ and $M$ is variable then, as mentioned
in the introduction, the problem is NP-hard. We now show that,
for fixed $d$, Theorem \ref{CCO} implies at once that the problem
is strongly polynomial time solvable (see also \cite{AFLS}).

\bc{Quadratic}
For every fixed $d$ there is a strongly polynomial time algorithm
that given $W\in\Z^{d\times n}$, encoded as $[n;\l W\r]$,
finds $x^*\in\{0,1\}^n$ maximizing the form $x^TW^TWx$.
\ec
\bpr
Let $S:=\{0,1\}^n$ and let $E:=\{{\bf 1}_1,\dots,{\bf 1}_n\}$ be the
set of unit vectors in $\R^n$. Then $P:=\conv(S)$ is just the $n$-cube
$[0,1]^n$ and hence $E$ covers all edge-directions of $P$.
A membership oracle for $S$ is easily and efficiently realizable and $x:=0\in S$
is an initial point. Also, $|E|$ and $\l E\r$ are polynomial in $n$,
and $E$ is easily and efficiently computable.

Now, for $i=1,\dots,d$ define $w_i\in\Z^n$ to be the $i$-th row of
the matrix $W$, that is, $w_{i,j}:=W_{i,j}$ for all $i,j$. Finally,
let $c:\R^d\longrightarrow\R$ be the squared $l_2$ norm given by
$c(y):=\|y\|_2^2:=\sum_{i=1}^d y_i^2$, and note that the comparison
of $c(y)$ and $c(z)$ can be done for $y,z\in\Z^d$ in time polynomial
in $\l y,z\r$ using a constant number of arithmetic operations, providing
a strongly polynomial time realization of a comparison oracle for $c$.

This translates the given quadratic programming problem into a convex
combinatorial optimization problem over $S$, which can be solved in
strongly polynomial time by applying the algorithm of
Theorem \ref{CCO} to $S$, $x=0$, $w_1,\dots,w_d$, $E$, and $c$.
\epr

\subsubsection{Matroids and Maximum Norm Spanning Trees}
\label{mamnst}

Optimization problems over matroids form a fundamental class of combinatorial
optimization problems. Here we discuss matroid bases, but everything works
for independent sets as well. Recall that a family $\cal B$ of subsets of
$\{1,\dots,n\}$ is the family of {\em bases} of a {\em matroid} if all members
of $\cal B$ have the same cardinality, called the {\em rank} of the matroid,
and for every $B,B'\in\cal B$ and $i\in B\setminus B'$ there is a $j\in B'$
such that $B\setminus\{i\}\cup\{j\}\in{\cal B}$. Useful models include the
{\em graphic matroid} of a graph $G$ with edge set $\{1,\dots,n\}$ and $\cal B$
the family of spanning forests of $G$, and the {\em linear matroid} of an
$m\times n$ matrix $A$ with $\cal B$ the family of sets of indices
of maximal linearly independent subsets of columns of $A$.

It is well known that linear combinatorial optimization over matroids
can be solved by the fast greedy algorithm \cite{Edm}. We now show that,
as a consequence of Theorem \ref{CCO}, convex combinatorial optimization
over a matroid presented by a membership oracle can be solved in strongly
polynomial time as well (see also \cite{HT,Onn1}). We state the result for
bases, but the analogous statement for independent sets hold as well.
We say that $S\subseteq\{0,1\}^n$ is the {\em set of bases of a matroid}
if it is the set of indicators of the family $\cal B$ of bases of some
matroid, in which case we call $\conv(S)$ the {\em matroid base polytope}.

\bc{Matroids}
For every fixed $d$ there is a strongly polynomial time algorithm that, given
set $S\subseteq\{0,1\}^n$ of bases of a matroid presented by a membership
oracle, $x\in S$, $w_1,\dots,w_d\in\Z^n$, and convex functional
$c:\R^d\longrightarrow\R$ presented by a comparison oracle, encoded as
$[n;\l x,w_1,\dots,w_d\r]$, solves the convex matroid optimization problem
$$\max\, \{c(w_1 z,\dots,w_d z):\ z\in S\}\ .$$
\ec
\bpr
Let $E:=\{{\bf 1}_i-{\bf 1}_j: 1\leq i< j\leq n\}$ be the set of
differences of pairs of unit vectors in $\R^n$. We claim that $E$
covers all edge-directions of the matroid base polytope $P:=\conv(S)$.
Consider any edge $e=[y,y']$ of $P$ with $y,y'\in S$ and let
$B:=\supp(y)$ and $B':=\supp(y')$ be the corresponding bases.
Let $h\in\R^n$ be a linear functional uniquely maximized over $P$ at $e$.
If $B\setminus B'=\{i\}$ is a singleton then $B'\setminus B=\{j\}$ is a
singleton as well in which case $y-y'={\bf 1}_i-{\bf 1}_j$ and we are done.
Suppose then, indirectly, that it is not, and pick an element $i$ in the
symmetric difference $B\Delta B':=(B\setminus B')\cup(B'\setminus B)$ of
minimum value $h_i$. Without loss of generality assume $i\in B\setminus B'$.
Then there is a $j\in B'\setminus B$ such that $B'':=B\setminus\{i\}\cup\{j\}$
is also a basis. Let $y''\in S$ be the indicator of $B''$. Now $|B\Delta B'|>2$
implies that $B''$ is neither $B$ nor $B'$. By the choice of $i$ we have
$hy''=hy-h_i+h_j\geq hy$. So $y''$ is also a maximizer of $h$ over
$P$ and hence $y''\in e$. But no $\{0,1\}$-vector is a convex
combination of others, a contradiction.

Now, $|E|={n\choose 2}$ and $E\subset \{-1,0,1\}^n$ imply
that $|E|$ and $\l E\r$ are polynomial in $n$. Moreover, $E$ can
be easily computed in strongly polynomial time. Therefore,
applying the algorithm of Theorem \ref{CCO} to the
given data and the set $E$, the convex discrete optimization
problem over $S$ can be solved in strongly polynomial time.
\epr

One important application of Corollary \ref{Matroids} is a polynomial time
algorithm for computing the {\em universal Gr\"obner basis} of any system
of polynomials having a finite set of common zeros in fixed (but arbitrary)
number of variables, as well as the construction of the {\em state polyhedron}
of any member of the {\em Hilbert scheme}, see \cite{BOT,OSt}.
Other important applications are in the field of {\em algebraic statistics}
\cite{PRW}, in particular for {\em optimal experimental design}.
These applications are beyond our scope here and will be discussed elsewhere.

\vskip.2cm
Here is another concrete example of a convex matroid optimization application.

\be{Matroid1}{\sc (maximum norm spanning tree).}
Fix any positive integer $d$. Let $\|\cdot\|_p:\R^d\longrightarrow\R$ be
the $l_p$~norm given by $\|x\|_p:=(\sum_{i=1}^d |x_i|^p)^{1\over p}$
for $1\leq p<\infty$ and $\|x\|_{\infty}:=\max_{i=1}^d |x_i|$.
Let $G$ be a connected graph with edge set $N:=\{1,\dots,n\}$.
For $j=1,\dots,n$ let $u_j\in\Z^d$ be a weight vector representing
the values of edge $j$ under some $d$ criteria. The weight
of a subset $T\subseteq N$ is the sum $\sum_{j\in T}u_j$ representing
the total values of $T$ under the $d$ criteria. The problem is to find a
spanning tree $T$ of $G$ whose weight has maximum $l_p$~norm, that is,
a spanning tree $T$ maximizing $\|\sum_{j\in T}u_j\|_p$.

Define $w_1,\dots,w_d\in\Z^n$ by $w_{i,j}:=u_{j,i}$ for $i=1,\dots,d$,
$j=1,\dots,n$. Let $S\subseteq\{0,1\}^n$ be the set of indicators of
spanning trees of $G$. Then, in time polynomial in $n$, a membership
oracle for $S$ is realizable, and an initial $x\in S$ is obtainable
as the indicator of any greedily constructible spanning tree $T$.
Finally, define the convex functional $c:=\|\cdot\|_p$. Then for
most common values $p=1,2,\infty$, and in fact for any $p\in\N$,
the comparison of $c(y)$ and $c(z)$ can be done for $y,z\in\Z^d$ in
time polynomial in $\l y,z,p\r$ by computing and comparing the integer
valued $p$-th powers $\|y\|^p_p$ and $\|z\|^p_p$. Thus, by
Corollary \ref{Matroids}, this problem is solvable in time polynomial
in $\l u_1,\dots,u_n,p \r$.
\ee

\newpage
\section{Linear N-fold Integer Programming}
\label{lnfip}

In this section we develop a theory of linear
{\em $n$-fold integer programming}, which leads to the polynomial
time solution of broad classes of linear integer programming
problems in variable dimension. This will be extended to
convex $n$-fold integer programming in \S \ref{cip}.

In \S \ref{oaalo} we describe an adaptation of a result
of \cite{SW} involving an oriented version of the augmentation
oracle of \S \ref{fmtlo}. In \S \ref{gbalip} we discuss
Graver bases and their application to linear integer programming.
In \S \ref{gbonfm} we show that Graver bases of $n$-fold matrices can be
computed efficiently. In \S \ref{lnfipipt} we combine the preparatory
statements from \S \ref{oaalo}, \S \ref{gbalip}, and \S \ref{gbonfm},
and prove the main result of this section, asserting that
linear $n$-fold integer programming is polynomial time solvable.
We conclude with some applications in \S \ref{salip}.

Here and in \S \ref{cip} we concentrate on discrete optimization
problems over a set $S$ presented as the set of integer
points satisfying an explicitly given system of linear inequalities.
Without loss of generality we may and will assume that $S$ is given either
in standard form $S:=\{x\in\N^n:\ Ax=b\}$ where $A\in\Z^{m\times n}$
and $b\in\Z^m$, or in the form
$$S\quad:=\quad\{x\in\Z^n:\ Ax=b,\ l\leq x\leq u\}$$
where $l,u\in\Z_{\infty}^n$ and $\Z_{\infty}=\Z\uplus\{\pm\infty\}$,
where some of the variables are bounded below or above and some
are unbounded. Thus, $S$ is no longer presented by an oracle,
but by the explicit data $A,b$ and possibly $l,u$.
In this setup we refer to discrete optimization over $S$ also as
{\em integer programming} over $S$. As usual, an algorithm solving the
problem must either provide an $x\in S$ maximizing $wx$ over $S$,
or assert that none exists (either because $S$ is empty or because
the objective function is unbounded over the underlying polyhedron).
We will sometimes assume that an initial point $x\in S$ is given, in which
case $b$ will be computed as $b:=Ax$ and not be part of the input.

\subsection{Oriented Augmentation and Linear Optimization}
\label{oaalo}

We have seen in \S \ref{fmtlo} that an augmentation oracle presentation
of a finite set $S\subset\Z^n$ enables to solve the linear discrete
optimization problem over $S$. However, the running time of the
algorithm of Lemma \ref{Augmentation} which demonstrated
this, was polynomial in the unary length of the radius $\rho(S)$
of the feasible set rather than in its binary length.

In this subsection we discuss a recent result of \cite{SW} and show that,
when $S$ is presented by a suitable stronger oriented version of the
augmentation oracle, the linear optimization problem can be solved by
a much faster algorithm, whose running time is in fact polynomial in the
binary length $\l\rho(S)\r$. The key idea behind this algorithm is that
it gives preference to augmentations along interior points of
$\conv(S)$ staying far off its boundary. It is inspired by and
extends the combinatorial interior point algorithm of \cite{WZ}.

For any vector $g\in\R^n$, let $g^+,g^-\in\R_+^n$ denote its
{\em positive} and {\em negative} parts, defined by
$g^+_j:=\max\{g_j,0\}$ and $g^-_j:=-\min\{g_j,0\}$ for
$j=1,\dots,n$. Note that both $g^+,g^-$ are nonnegative,
$\supp(g)=\supp(g^+)\biguplus \supp(g^-)$, and $g=g^+-g^-$.

An {\em oriented augmentation oracle} for a set $S\subset\Z^n$ is one that,
queried on $x\in S$ and $w_+,w_-\in\Z^n$, either returns an
{\em augmenting vector} $g\in\Z^n$, defined to be one satisfying
$x+g\in S$ and $w_+g^+-w_-g^->0$, or asserts that none exists.

Note that this oracle involves {\em two} linear functionals $w_+,w_-\in\Z^n$
rather than one ($w_+,w_-$ are two distinct independent vectors and
{\em not} the positive and negative parts of one vector).
The conditions on an augmenting vector $g$ indicate that it is
a feasible direction and has positive value under the nonlinear objective
function determined by $w_+,w_-$. Note that this oracle is indeed
stronger than the augmentation oracle of \S \ref{fmtlo}:
to answer a query $x\in S$, $w\in\Z^n$ to the latter, set $w_+:=w_-:=w$,
thereby obtaining $w_+g^+-w_-g^-=wg$ for all $g$, and query
the former on $x,w_+,w_-$; if it replies with an augmenting vector
$g$ then reply with the better point ${\hat x}:=x+g$, whereas if it
asserts that no $g$ exists then assert that no better point exists.

The following lemma is an adaptation of the result of \cite{SW} concerning
sets of the form $S:=\{x\in\Z^n:Ax=b,\ 0\leq x\leq u\}$ of nonnegative integer points
satisfying equations and upper bounds. However, the pair $A,b$ is neither
explicitly needed nor does it affect the running time of the algorithm
underlying the lemma. It suffices that $S$ is of that form. Moreover,
an arbitrary lower bound vector $l$ rather than $0$ can be included.
So it suffices to assume that $S$ coincides with the intersection
of its affine hull and the set of integer points in a box, that is,
$S=\aff(S)\cap \{x\in\Z^n:l\leq x\leq u\}$ where $l,u\in\Z^n$.
We now describe and prove the algorithm of \cite{SW} adjusted to any
lower and upper bounds $l,u$.

\bl{OrientedAugmentation}
There is a polynomial time algorithm that, given vectors $l,u\in\Z^n$, set\break
$S\subset\Z^n$ satisfying $S=\aff(S)\cap\{z\in\Z^n:l\leq z\leq u\}$
and presented by an oriented augmentation oracle, $x\in S$,
and $w\in\Z^n$, encoded as $[\l l,u,x,w \r]$, provides an optimal solution
$x^*\in S$ to the linear discrete optimization problem $\max\{wz:z\in S\}$.
\el
\bpr
We start with some strengthening adjustments to the oriented augmentation
oracle. Let $\rho:=\max\{\|l\|_{\infty},\|u\|_{\infty}\}$ be an upper bound
on the radius of $S$. Then any augmenting vector $g$ obtained from the
oriented augmentation oracle when queried on $y\in S$ and $w_+,w_-\in\Z^n$,
can be made in polynomial time to be {\em exhaustive}, that is, to satisfy\break
$y+2g\notin S$ (which means that no longer augmenting step in direction $g$
can be taken). Indeed, using binary search, find the largest
$r\in\{1,\dots,2\rho\}$ for which $l\leq y+rg\leq u$; then
$S=\aff(S)\cap\{z\in\Z^n:l\leq z\leq u\}$ implies $y+rg\in S$ and hence
we can replace $g:=rg$. So from here on we will assume that if there
is an augmenting vector then the oracle returns an exhaustive one.
Second, let $\R_{\infty}:=\R\uplus\{\pm\infty\}$ and for any vector
$v\in\R^n$ let $v^{-1}\in\R_{\infty}^n$ denote its entry-wise
reciprocal defined by $v^{-1}_i:={1\over v_i}$ if $v_i\neq 0$ and
$v^{-1}_i:=\infty$ if $v_i=0$. For any $y\in S$, the vectors
$(y-l)^{-1}$ and $(u-y)^{-1}$ are the reciprocals of the
``entry-wise distance" of $y$ from the given lower and upper bounds.
The algorithm will query the oracle on triples $y,w_+,w_-$
with $w_+:=w-\mu(u-y)^{-1}$ and $w_-:=w+\mu(y-l)^{-1}$ where $\mu$
is a suitable positive scalar and $w$ is the input linear functional. The fact
that such $w_+,w_-$ may have infinite entries does not cause any problem:
indeed, if $g$ is an augmenting vector then $y+g\in S$ implies that
$g^+_i=0$ whenever $y_i=u_i$ and $g^-_i=0$ whenever $l_i=y_i$,
so each infinite entry in $w_+$ or $w_-$ occurring in the
expression $w_+g^+-w_-g^-$ is multiplied by $0$ and hence zeroed out.

The algorithm proceeds in phases. Each phase $i$ starts with a feasible
point $y_{i-1}\in S$ and performs repeated augmentations using the oriented
augmentation oracle, terminating with a new feasible point $y_i\in S$
when no further augmentations are possible. The queries to the oracle
make use of a positive scalar parameters $\mu_i$ fixed throughout the
phase. The first phase ($i$=1) starts with the input point $y_0:=x$ and sets
$\mu_1:=\rho\,\|w\|_{\infty}$. Each further phase $i\geq 2$ starts with the
point $y_{i-1}$ obtained from the previous phase and sets the parameter
value $\mu_i:={1\over 2}\mu_{i-1}$ to be half its value in the previous phase.
The algorithm terminates at the end of the first phase $i$ for which
$\mu_i<{1\over n}$, and outputs $x^*:=y_i$. Thus, the number of phases
is at most $\lceil\log_2(2n\rho\|w\|_{\infty})\rceil$
and hence polynomial in $\l l,u,w \r$.

We now describe the $i$-th phase which determines $y_i$ from $y_{i-1}$.
Set $\mu_i:={1\over 2}\mu_{i-1}$ and ${\hat y}:=y_{i-1}$.
Iterate the following: query the strengthened
oriented augmentation oracle on ${\hat y}$, $w_+:=w-\mu_i(u-{\hat y})^{-1}$,
and $w_-:=w+\mu_i({\hat y}-l)^{-1}$; if the oracle returns an exhaustive
augmenting vector $g$ then set ${\hat y}:={\hat y}+g$ and repeat, whereas
if it asserts that there is no augmenting vector then set $y_i:={\hat y}$
and complete the phase. If $\mu_i\geq {1\over n}$ then
proceed to the $(i+1)$-th phase, else
output $x^*:=y_i$ and terminate the algorithm.

It remains to show that the output of the algorithm is indeed an optimal
solution and that the number of iterations (and hence calls to the oracle)
in each phase is polynomial in the input. For this we need the following
facts, the easy proofs of which are omitted:
\begin{enumerate}
\item
For every feasible $y\in S$ and direction $g$ with
$y+g\in S$ also feasible, we have
$$(u-y)^{-1}g^+ +(y-l)^{-1}g^-\ \ \leq\ \ n\ \ .$$
\item
For every $y\in S$ and direction $g$ with
$y+g\in S$ but $y+2g\notin S$, we have
$$(u-y)^{-1}g^+ +(y-l)^{-1}g^-\ \ >\ \ {1\over2}\ \  .$$
\item
For every feasible $y\in S$, direction $g$ with $y+g\in S$ also feasible, and
$\mu>0$, setting $w_+:=w-\mu(u-y)^{-1}$ and $w_-:=w+\mu(y-l)^{-1}$ we have
$$w_+g^+-w_-g^-\ \ =\ \ wg-\mu\left((u-y)^{-1}g^+ +(y-l)^{-1}g^-\right)\ \ .$$
\end{enumerate}
Now, consider the last phase $i$ with $\mu_i<{1\over n}$, let
$x^*:=y_i:={\hat y}$ be the output of the algorithm at the end
of this phase, and let $\hat x\in S$ be any optimal solution.
Now, the phase is completed when the oracle, queried on the triple ${\hat y}$,
$w_+=w-\mu_i(u-{\hat y})^{-1}$, and $w_-=w+\mu_i({\hat y}-l)^{-1}$,
asserts that there is no augmenting vector. In particular, setting
$g:={\hat x}-{\hat y}$, we find $w_+g^+-w_-g^-\leq 0$ and hence,
by facts 1 and 3 above,
$$w{\hat x}-wx^*\ =\ wg
\ \leq \ \mu_i \left((u-{\hat y})^{-1}g^+ +({\hat y}-l)^{-1}g^-\right)
\ <\ {1\over n}\, n\ =\ 1\ .$$
Since $w{\hat x}$ and $wx^*$ are integer, this implies that in fact
$w{\hat x}-wx^*\leq 0$ and hence the output $x^*$ of the algorithm
is indeed an optimal solution to the given optimization problem.

Next we bound the number of iterations in each phase $i$ starting
from $y_{i-1}\in S$. Let again $\hat x\in S$ be any optimal solution.
Consider any iteration in that phase, where the oracle is queried on
${\hat y}$, $w_+=w-\mu_i(u-{\hat y})^{-1}$, and $w_-=w+\mu_i({\hat y}-l)^{-1}$,
and returns an exhaustive augmenting vector $g$. We will now show that
\begin{equation}\label{e1}
w({\hat y}+g)-w{\hat y}\quad \geq \quad{1\over 4n}(w{\hat x}-wy_{i-1})\ \ ,
\end{equation}
that is, the increment in the objective value from ${\hat y}$ to the
augmented point ${\hat y}+g$ is at least ${1\over 4n}$ times
the difference between the optimal objective value $w{\hat x}$ and the
objective value $wy_{i-1}$ of the point $y_{i-1}$  at the beginning of phase $i$.
This shows that at most $4n$ such increments (and hence iterations)
can occur in the phase before it is completed.

To establish (\ref{e1}), we show that $wg\geq {1\over2}\mu_i$ and
$w{\hat x}-wy_{i-1}\leq 2n\mu_i$. For the first inequality,
note that $g$ is an exhaustive augmenting vector and so $w_+g^+-w_-g^->0$
and ${\hat y}+2g\notin S$ and hence, by facts 2 and 3,
$wg>\mu_i((u-{\hat y})^{-1}g^+ +({\hat y}-l)^{-1}g^-)>{1\over2}\mu_i$.
We proceed with the second inequality. If $i=1$ (first phase) then this
indeed holds since $w{\hat x}-wy_0\leq 2n\rho\|w\|_{\infty} =2n\mu_1$.
If $i\geq 2$, let ${\tilde w}_+:=w-\mu_{i-1}(u-y_{i-1})^{-1}$ and
${\tilde w}_-:=w+\mu_{i-1}(y_{i-1}-l)^{-1}$. The $(i-1)$-th phase was
completed when the oracle, queried on the triple $y_{i-1}$, ${\tilde w}_+$,
and ${\tilde w}_-$, asserted that there is no augmenting vector.
In particular, for ${\tilde g}:={\hat x}-y_{i-1}$, we find
${\tilde w}_+{\tilde g}^+ - {\tilde w}_-{\tilde g}^-\leq 0$ and so,
by facts 1 and 3,
$$w{\hat x}-wy_{i-1}=w{\tilde g}\leq\mu_{i-1}
\left((u-y_{i-1})^{-1}{\tilde g}^+ +(y_{i-1}-l)^{-1}{\tilde g}^-)\right)
\leq\mu_{i-1}n=2n\mu_i\,.\ \, \square$$

\subsection{Graver Bases and Linear Integer Programming}
\label{gbalip}

We now come to the definition of a fundamental object introduced by
Graver in \cite{Gra}. The {\em Graver basis} of an integer matrix $A$
is a canonical finite set $\G(A)$ that can be defined
as follows. Define a partial order $\sqsubseteq$ on $\Z^n$
which extends the coordinate-wise order $\leq$ on $\N^n$ as follows: for two
vectors $u,v\in\Z^n$ put $u\sqsubseteq v$ and say that $u$ is {\em conformal}
to $v$ if $|u_i|\leq |v_i|$ and $u_iv_i\geq 0$ for $i=1,\ldots,n$,
that is, $u$ and $v$ lie in the same orthant of $\R^n$ and each component of
$u$ is bounded by the corresponding component of $v$ in absolute value.
It is not hard to see that $\sqsubseteq$ is a {\em well} partial ordering
(this is basically Dickson's lemma) and hence every subset of $\Z^n$ has
finitely-many $\sqsubseteq$-minimal elements.
Let $\L(A):=\{x\in\Z^n:\ Ax=0\}$ be the lattice of linear
integer dependencies on $A$. The {\em Graver basis} of $A$
is defined to be the set $\G(A)$ of all $\sqsubseteq$-minimal
vectors in $\L(A)\setminus\{0\}$.

Note that if $A$ is an $m\times n$ matrix then its Graver basis consist of
vectors in $\Z^n$. We sometimes write $\G(A)$ as a suitable $|\G(A)|\times n$
matrix whose rows are the Graver basis elements. The Graver basis is centrally
symmetric ($g\in\G(A)$ implies $-g\in\G(A)$); thus, when listing a Graver
basis we will typically give one of each antipodal pair and prefix the set
(or matrix) by $\pm$. Any element of the Graver basis is primitive
(its entries are relatively prime integers). Every circuit of $A$
(nonzero primitive minimal support element of $\L(A)$) is in $\G(A)$; in fact,
if $A$ is totally unimodular then $\G(A)$ coincides with the set
of circuits (see \S \ref{cipotus} in the sequel for more details on this).
However, in general $\G(A)$ is much larger. For more details on Graver bases
and their connection to Gr\"obner bases see Sturmfels \cite{Stu} and for the
currently fastest procedure for computing them see \cite{Hem,HHM}.

Here is a quick simple example; we will see more structured and
complex examples later on. Consider the $1\times 3$ matrix
$A:=(1,2,1)$. Then its Graver basis can be shown to be the set
$\G(A)=\pm\{(2,-1,0),(0,-1,2),(1,0,-1),(1,-1,1)\}$.
The first three elements (and their antipodes) are the circuits of $A$;
already in this small example non-circuits appear as well: the fourth element
(and its antipode) is a primitive linear integer dependency whose support
is not minimal.

We now show that when we do have access to the Graver basis,
it can be used to solve linear integer programming. We will extend
this in \S \ref{cip}, where we show that the
Graver basis enables to solve convex integer programming as well.
In \S \ref{gbonfm} we will show that there are important
classes of matrices for which the Graver basis is indeed accessible.

First, we need a simple property of Graver bases. A finite sum
$u:=\sum_i v_i$ of vectors $v_i\in\R^n$ is {\em conformal} if each
summand is conformal to the sum, that is, $v_i\sqsubseteq u$ for all $i$.

\bl{Conformal}
Let $A$ be any integer matrix. Then any $h\in\L(A)\setminus\{0\}$
can be written as a conformal sum $h:=\sum g_i$ of
(not necessarily distinct) Graver basis elements $g_i\in\G(A)$.
\el
\bpr
By induction on the well partial order $\sqsubseteq$.
Recall that $\G(A)$ is the set of $\sqsubseteq$-minimal elements
in $\L(A)\setminus\{0\}$. Consider any $h\in\L(A)\setminus\{0\}$.
If it is $\sqsubseteq$-minimal then $h\in\G(A)$ and we are done.
Otherwise, there is a $h'\in\G(A)$ such that $h'\sqsubset h$.
Set $h'':=h-h'$. Then $h''\in\L(A)\setminus\{0\}$ and $h''\sqsubset h$,
so by induction there is a conformal sum $h''=\sum_i g_i$ with $g_i\in\G(A)$
for all $i$. Now $h=h'+\sum_i g_i$ is the desired conformal sum of $h$.
\epr

The next lemma shows the usefulness of Graver bases for oriented augmentation.
\bl{GraverOriented}
Let $A$ be an $m\times n$ integer matrix with Graver basis $\G(A)$
and let $l,u\in\Z_{\infty}^n$, $w_+,w_-\in\Z^n$, and $b\in\Z^m$.
Suppose $x\in T:=\{y\in\Z^n:Ay=b,l\leq y\leq u\}$.
Then for every $g\in\Z^n$ which satisfies $x+g\in T$ and $w_+g^+-w_-g^->0$ there
exists an  element ${\hat g}\in\G(A)$ with ${\hat g}\sqsubseteq g$
which also satisfies $x+{\hat g}\in T$ and $w_+{\hat g}^+-w_-{\hat g}^->0$.
\el
\bpr
Suppose $g\in\Z^n$ satisfies the requirements. Then $Ag=A(x+g)-Ax=b-b=0$
since $x,x+g\in T$. Thus, $g\in\L(A)\setminus\{0\}$
and hence, by Lemma \ref{Conformal}, there is
a conformal sum $g=\sum_i h_i$ with $h_i\in\G(A)$ for all $i$.
Now, $h_i\sqsubseteq g$ is equivalent to $h_i^+\leq g^+$ and
$h_i^-\leq g^-$, so the conformal sum $g=\sum_i h_i$ gives
corresponding sums of the positive and negative parts
$g^+=\sum_i h_i^+$ and $g^-=\sum_i h_i^-$. Therefore we obtain
$$
0\ <\ w_+g^+ - w_-g^- \ =\  w_+\sum_i h_i^+ - w_-\sum_i h_i^-
\ = \ \sum_i (w_+h_i^+ - w_-h_i^-)
$$
which implies that there is some $h_i$ in this sum with
$w_+h_i^+ - w_-h_i^->0$. Now, $h_i\in\G(A)$ implies
$A(x+h_i)=Ax=b$. Also, $l\leq x,x+g \leq u$ and $h_i\sqsubseteq g$
imply that $l\leq x+h_i\leq u$. So $x+h_i\in T$.
Therefore the vector ${\hat g}:=h_i$ satisfies the claim.
\epr

We can now show that the Graver basis enables to solve
linear integer programming in polynomial time
provided an initial feasible point is available.

\bt{GraverIP}
There is a polynomial time algorithm that, given $A\in\Z^{m\times n}$,
its Graver basis $\G(A)$, $l,u\in\Z_{\infty}^n$, $x,w\in\Z^n$ with $l\leq x\leq u$,
encoded as $[\l A,\G(A),l,u,x,w\r]$, solves the linear integer program
$\max\{wz:z\in\Z^n,Az=b,l\leq z\leq u\}$ with $b:=Ax$.
\et
\bpr
First, note that the objective function of the integer program is
unbounded if and only if the objective function of its relaxation
$\max\{wy:y\in\R^n,Ay=b,l\leq y\leq u\}$ is unbounded, which can be checked
in polynomial time using linear programming. If it is unbounded
then assert that there is no optimal solution and terminate the algorithm.

Assume then that the objective is bounded. Then, since the
program is feasible, it has an optimal solution. Furthermore,
(as basically follows from Cramer's rule, see e.g. \cite[Theorem 17.1]{Sch})
it has an optimal $x^*$ satisfying $|x^*_j|\leq \rho$ for all $j$,
where $\rho$ is an easily computable integer upper bound whose binary length
$\l\rho\r$ is polynomially bounded in $\l A,l,u,x \r$.
For instance, $\rho:=(n+1)(n+1)!r^{n+1}$ will do, with $r$ the maximum among
$\max_i|\sum_jA_{i,j}x_j|$, $\max_{i,j}|A_{i,j}|$,
$\max\{|l_j|:|l_j|<\infty\}$, and $\max\{|u_j|: |u_j|<\infty\}$.

Let $T:=\{y\in\Z^n:Ay=b,\ l\leq y\leq u\}$ and $S:=T\cap[-\rho,\rho]^n$.
Then our linear integer programming problem now reduces
to linear discrete optimization over $S$.
Now, an oriented augmentation oracle for $S$ can be simulated
in polynomial time using the given Graver basis $\G(A)$ as follows:
given a query $y\in S$ and $w_+,w_-\in\Z^n$, search for $g\in\G(A)$
which satisfies $w_+g^+-w_-g^->0$ and $y+g\in S$; if there is such
a $g$ then return it as an augmenting vector, whereas if there is
no such $g$ then assert that no augmenting vector exists.
Clearly, if this simulated oracle returns a vector $g$ then it is
an augmenting vector. On the other hand, if there exists an augmenting
vector $g$ then $y+g\in S\subseteq T$ and
$w_+g^+-w_-g^->0$ imply by Lemma \ref{GraverOriented} that there
is also a ${\hat g}\in\G(A)$ with ${\hat g}\sqsubseteq g$ such that
$w_+{\hat g}^+-w_-{\hat g}^->0$  and $y+{\hat g}\in T$. Since
$y,y+g\in S$ and ${\hat g}\sqsubseteq g$, we find that $y+{\hat g}\in S$
as well. Therefore the Graver basis contains an augmenting vector
and hence the simulated oracle will find and output one.

Define ${\hat l},{\hat u}\in\Z^n$ by
${\hat l}_j:=\max(l_j,-\rho),{\hat u}_j:=\min(u_j,\rho)$, $j=1,\dots,n$.
Then it is easy to see that
$S=\aff(S)\cap\{y\in\Z^n:{\hat l}\leq y\leq{\hat u}\}$.
Now apply the algorithm of Lemma \ref{OrientedAugmentation}
to ${\hat l},{\hat u}$, $S$, $x$, and $w$, using the above simulated
oriented augmentation oracle for $S$, and obtain in polynomial time a
vector $x^*\in S$ which is optimal to the linear discrete optimization
problem over $S$ and hence to the given linear integer program.
\epr

As a special case of Theorem \ref{GraverIP} we recover the following result
of \cite{DHOW} concerning linear integer programming in standard form
when the Graver basis is available.

\bt{SpecialGraverIP}
There is a polynomial time algorithm that, given matrix $A\in\Z^{m\times n}$,
its Graver basis $\G(A)$, $x\in\N^n$, and $w\in\Z^n$, encoded as
$[\l A,\G(A),x,w\r]$, solves the linear integer programming problem
$\max\{wz:z\in\N^n,\ Az=b\}$ where $b:=Ax$.
\et

\subsection{Graver Bases of N-fold Matrices}
\label{gbonfm}

As mentioned above, the Graver basis $\G(A)$ of an integer matrix $A$
contains all circuits of $A$ and typically many more elements.
While the number of circuits is already typically
exponential and can be as large as $n\choose{m+1}$, the number of Graver
basis elements is usually even larger and depends also on the entries
of $A$ and not only on its dimensions $m,n$. So unfortunately
it is typically very hard to compute $\G(A)$. However, we now show
that for the important and useful broad class of $n$-fold matrices,
the Graver basis is better behaved and can be computed in polynomial time.
Recall the following definition from the introduction.
Given an $(r+s)\times t$ matrix $A$, let $A_1$ be its $r\times t$ sub-matrix
consisting of the first $r$ rows and let $A_2$ be its $s\times t$
sub-matrix consisting of the last $s$ rows. We refer to $A$ explicitly
as {\em $(r+s)\times t$ matrix}, since the definition below depends also
on $r$ and $s$ and not only on the entries of $A$. The {\em $n$-fold matrix}
of an $(r+s)\times t$ matrix $A$ is then defined to be
the following $(r+ns)\times nt$ matrix,
$$A^{(n)}\quad:=\quad ({\bf 1}_n\otimes A_1)\oplus(I_n \otimes A_2)\quad=\quad
\left(
\begin{array}{ccccc}
  A_1    & A_1    & A_1    & \cdots & A_1    \\
  A_2  & 0      & 0      & \cdots & 0      \\
  0  & A_2      & 0      & \cdots & 0      \\
  \vdots & \vdots & \vdots & \ddots & \vdots \\
  0  & 0      & 0      & \cdots & A_2      \\
\end{array}
\right)\quad .
$$

We now discuss a recent result of \cite{SS}, which originates in \cite{AT},
and its extension in \cite{HS}, on the stabilization
of Graver bases of $n$-fold matrices. Consider vectors $x=(x^1,\ldots,x^n)$
with $x^k\in\Z^t$ for $k=1,\dots,n$. The {\em type} of $x$ is the number
$|\{k\,:\,x^k\neq 0\}|$ of nonzero components $x^k\in\Z^t$ of $x$.
The {\em Graver complexity} of an $(r+s)\times t$ matrix,
denoted $c(A)$, is defined to be the smallest $c\in\N\uplus\{\infty\}$
such that for all $n$, the Graver basis of $A^{(n)}$ consists of
vectors of type at most $c(A)$. We provide the proof of the following result
of \cite{HS,SS} stating that the Graver complexity is always finite.

\bl{GraverComplexity}
The Graver complexity $c(A)$ of any $(r+s)\times t$ integer matrix $A$ is finite.
\el
\bpr
Call an element $x=(x^1,\ldots,x^n)$ in the Graver basis of some $A^{(n)}$
{\em pure} if\break $x^i\in\G(A_2)$ for all $i$. Note that the type
of a pure $x\in\G(A^{(n)})$ is $n$. First, we claim that if
there is an element of type $m$ in some $\G(A^{(l)})$ then for
some $n\geq m$ there is a pure element in $\G(A^{(n)})$, and so it
will suffice to bound the type of pure elements.
Suppose there is an element of type $m$ in some $\G(A^{(l)})$.
Then its restriction to its $m$ nonzero components
is an element $x=(x^1,\ldots,x^m)$ in $\G(A^{(m)})$.
Let $x^i=\sum_{j=1}^{k_i} g_{i,j}$ be a conformal decomposition
of $x^i$ with $g_{i,j}\in\G(A_2)$ for all $i,j$, and let
$n:=k_1+\cdots+k_m\geq m$. Then $g:=(g_{1,1},\dots,g_{m,k_m})$ is
in $\G(A^{(n)})$, else there would be $\hat g\sqsubset g$ in $\G(A^{(n)})$
in which case the nonzero $\hat x$ with
${\hat x}^i:=\sum_{j=1}^{k_i}{\hat g}_{i,j}$ for all $i$
would satisfy ${\hat x}\sqsubset x$ and
${\hat x}\in\L(A^{(m)})$, contradicting $x\in\G(A^{(m)})$.
Thus $g$ is a pure element of type $n\geq m$, proving the claim.

We proceed to bound the type of pure elements. Let $\G(A_2)=\{g_1,\dots,g_m\}$
be the Graver basis of $A_2$ and let $G_2$ be the $t\times m$ matrix whose
columns are the $g_i$. Suppose $x=(x^1,\ldots,x^n)\in\G(A^{(n)})$ is pure for
some $n$. Let $v\in\N^m$ be the vector with\break $v_i:=|\{k:x^k=g_i\}|$ counting
the number of $g_i$ components of $x$ for each $i$. Then
$\sum_{i=1}^m v_i$ is equal to the type $n$ of $x$. Next, note that
$A_1G_2v=A_1(\sum_{k=1}^n x^k)=0$ and hence $v\in\L(A_1G_2)$. We claim that,
moreover, $v\in\G(A_1G_2)$. Suppose indirectly not. Then there is
${\hat v}\in\G(A_1G_2)$ with ${\hat v}\sqsubset v$, and it is easy to
obtain a nonzero ${\hat x}\sqsubset x$ from $x$ by zeroing out some components
so that ${\hat v}_i=|\{k:{\hat x}^k=g_i\}|$
for all $i$. Then $A_1(\sum_{k=1}^n {\hat x}^k)=A_1G_2{\hat v}=0$
and hence ${\hat x}\in\L(A^{(n)})$, contradicting $x\in\G(A^{(n)})$.

So the type of any pure element, and hence the Graver
complexity of $A$, is at most the largest value $\sum_{i=1}^m v_i$
of any nonnegative element $v$ of the Graver basis $\G(A_1G_2)$.
\epr

Using Lemma \ref{GraverComplexity} we now show how
to compute $\G(A^{(n)})$ in polynomial time.

\bt{GraverComputation}
For every fixed $(r+s)\times t$ integer matrix $A$ there is a strongly
polynomial time algorithm that, given $n\in\N$, encoded as $[n;n]$, computes
the Graver basis $\G(A^{(n)})$ of the $n$-fold matrix $A^{(n)}$. In particular,
the cardinality $|\G(A^{(n)})|$ and binary length\break $\l \G(A^{(n)})\r$
of the Graver basis of the $n$-fold matrix are polynomially bounded in $n$.
\et
\bpr
Let $c:=c(A)$ be the Graver complexity of $A$ and consider any $n\geq c$.
We show that the Graver basis of $A^{(n)}$ is the union of
$n\choose c$ suitably embedded copies of the Graver basis of $A^{(c)}$.
For every $c$ indices $1\leq k_1<\dots <k_c\leq n$ define a map
$\phi_{k_1,\dots,k_c}$ from $\Z^{ct}$ to $\Z^{nt}$ sending
$x=(x^1,\ldots,x^c)$ to $y=(y^1,\ldots,y^n)$ with
$y^{k_i}:=x^i$ for $i=1,\dots,c$ and $y^k:=0$ for $k\notin\{k_1,\dots,k_c\}$.
We claim that $\G(A^{(n)})$ is the union of the images
of $\G(A^{(c)})$ under the $n\choose c$ maps $\phi_{k_1,\dots,k_c}$
for all $1\leq k_1<\dots <k_c\leq n$, that is,
\begin{equation}\label{e2}
\G(A^{(n)})\quad=\quad \bigcup_{1\leq k_1<\dots <k_c\leq n}
\phi_{k_1,\dots,k_c}(\G(A^{(c)}))\quad .
\end{equation}
If $x=(x^1,\ldots,x^c)\in\G(A^{(c)})$ then $x$ is a
$\sqsubseteq$-minimal nonzero element of $\L(A^{(c)})$, implying that
$\phi_{k_1,\dots,k_c}(x)$ is a $\sqsubseteq$-minimal nonzero element
of $\L(A^{(n)})$ and therefore we have $\phi_{k_1,\dots,k_c}(x)\in\G(A^{(n)})$.
So the right-hand side of (\ref{e2}) is contained in the
left-hand side. Conversely, consider any $y\in\G(A^{(n)})$. Then,
by Lemma \ref{GraverComplexity}, the type of $y$ is at most $c$,
so there are indices $1\leq k_1<\dots <k_c\leq n$ such that all
nonzero components of $y$ are among those of the reduced vector
$x:=(y^{k_1},\ldots,y^{k_c})$ and therefore $y=\phi_{k_1,\dots,k_c}(x)$.
Now, $y\in\G(A^{(n)})$ implies that $y$ is a $\sqsubseteq$-minimal
nonzero element of $\L(A^{(n)})$ and hence $x$ is a $\sqsubseteq$-minimal
nonzero element of $\L(A^{(c)})$. Therefore $x\in\G(A^{(c)})$
and $y\in\phi_{k_1,\dots,k_g}(\G(A^{(c)}))$. So the
left-hand side of (\ref{e2}) is contained in the right-hand side.

Since $A$ is fixed we have that $c=c(A)$ and $\G(A^{(c)})$ are constant. Then (\ref{e2})
implies that $|\G(A^{(n)})|\leq{n\choose c}|\G(A^{(c)})|=O(n^c)$. Moreover,
every element of $\G(A^{(n)})$ is an $nt$-dimensional vector
$\phi_{k_1,\dots,k_c}(x)$ obtained by appending zero components to some
$x\in\G(A^{(c)})$ and hence has linear binary length $O(n)$. So the binary
length of the entire Graver basis $\G(A^{(n)})$ is $O(n^{c+1})$. Thus, the
${n\choose c}=O(n^c)$ images $\phi_{k_1,\dots,k_c}(\G(A^{(c)}))$ and their
union $\G(A^{(n)})$ can be computed in strongly polynomial time, as claimed.
\epr

\be{Example2}
Consider the $(2+1)\times 2$ matrix $A$ with $A_1:=I_2$ the $2\times 2$
identity and $A_2:=(1,1)$. Then $\G(A_2)=\pm(1,-1)$ and $\G(A_1G_2)=\pm(1,1)$
from which the Graver complexity of $A$ can be concluded to be $c(A)=2$
(see the proof of Lemma \ref{GraverComplexity}). The $2$-fold matrix of
$A$ and its Graver basis, consisting of two antipodal vectors only, are
$$
A^{(2)}
\ =\ \left(
\begin{array}{cccc}
  1 & 0 & 1 & 0 \\
  0 & 1 & 0 & 1 \\
  1 & 1 & 0 & 0 \\
  0 & 0 & 1 & 1 \\
\end{array}
\right)\,,\quad\quad
\G(A^{(2)})
\ =\ \pm \left(
\begin{array}{cccc}
  1 & -1 & -1 & 1 \\
\end{array}
\right)\quad.
$$
By Theorem \ref{GraverComputation}, the Graver basis of the $4$-fold matrix
$A^{(4)}$ is computed to be the union of the images of the $6={4\choose 2}$
maps $\phi_{k_1,k_2}:\Z^{2\cdot 2}\longrightarrow\Z^{4\cdot 2}$
for $1\leq k_1<k_2\leq 4$, getting
$$
A^{(4)}=\!\left(
\begin{array}{cccccccc}
  1 & 0 & 1 & 0 & 1 & 0 & 1 & 0 \\
  0 & 1 & 0 & 1 & 0 & 1 & 0 & 1 \\
  1 & 1 & 0 & 0 & 0 & 0 & 0 & 0 \\
  0 & 0 & 1 & 1 & 0 & 0 & 0 & 0 \\
  0 & 0 & 0 & 0 & 1 & 1 & 0 & 0 \\
  0 & 0 & 0 & 0 & 0 & 0 & 1 & 1 \\
\end{array}
\right)\!\!, \G(A^{(4)})
=\pm\!\left(
\begin{array}{cccccccc}
  1 & -1 & -1 & 1 & 0 & 0 & 0 & 0 \\
  1 & -1 & 0 & 0 & -1 & 1 & 0 & 0 \\
  1 & -1 & 0 & 0 & 0 & 0 & -1 & 1 \\
  0 & 0 & 1 & -1 & -1 & 1 & 0 & 0 \\
  0 & 0 & 1 & -1 & 0 & 0 & -1 & 1 \\
  0 & 0 & 0 & 0 & 1 & -1 & -1 & 1 \\
\end{array}
\right)\!.
$$
\ee

\subsection{Linear N-fold Integer Programming in Polynomial Time}
\label{lnfipipt}

We now proceed to provide a polynomial time algorithm for linear
integer programming over $n$-fold matrices. First, combining the results
of \S \ref{gbalip} and \S \ref{gbonfm}, we get at once the following
polynomial time algorithm for converting any feasible point to an optimal one.

\bl{NFoldAugmentation}
For every fixed $(r+s)\times t$ integer matrix $A$ there is a polynomial
time algorithm that, given $n\in\N$, $l,u\in\Z_{\infty}^{nt}$, $x,w\in\Z^{nt}$
satisfying $l\leq x\leq u$, encoded as $[\l l,u,x,w\r]$, solves the
linear $n$-fold integer programming problem with $b:=A^{(n)}x$,
$$\max\,\{wz\ :\ z\in\Z^{nt},\ A^{(n)}z=b,\ l\leq z\leq u\}\ .$$
\el
\bpr
First, apply the polynomial time algorithm of
Theorem \ref{GraverComputation} and compute the Graver basis
$\G(A^{(n)})$ of the $n$-fold matrix $A^{(n)}$. Then apply the
polynomial time algorithm of Theorem \ref{GraverIP}
to the data $A^{(n)}$, $\G(A^{(n)})$, $l,u,x$ and $w$.
\epr

Next we show that an initial feasible point
can also be found in polynomial time.

\bl{Feasible}
For every fixed $(r+s)\times t$ integer matrix $A$ there is a polynomial
time algorithm that, given $n\in\N$, $l,u\in\Z_{\infty}^{nt}$, and
$b\in\Z^{r+ns}$, encoded as $[\l l,u,b\r]$, either finds an $x\in\Z^{nt}$
satisfying $l\leq x\leq u$ and $A^{(n)}x=b$ or asserts that none exists.
\el
\bpr
If $l\not\leq u$ then assert that there is no feasible point and terminate
the algorithm. Assume then that $l\leq u$ and determine some $x\in\Z^{nt}$
with $l\leq x\leq u$ and $\l x\r\leq \l l,u\r$. Now, introduce $n(2r+2s)$
auxiliary variables to the given $n$-fold integer program and denote
by ${\hat x}$ the resulting vector of $n(t+2r+2s)$ variables. Suitably
extend the lower and upper bound vectors to ${\hat l},{\hat u}$ by setting
${\hat l}_j:=0$ and ${\hat u}_j:=\infty$ for each auxiliary variable
${\hat x}_j$. Consider the auxiliary integer program of finding an integer
vector $\hat x$ that minimizes the sum of auxiliary variables subject to
the lower and upper bounds ${\hat l}\leq {\hat x}\leq {\hat u}$ and the
following system of equations, with $I_r$ and $I_s$ the $r\times r$ and
$s\times s$ identity matrices,
$$\left(
\begin{array}{cccccccccccccccc}
A_1 & I_r & -I_r & 0 & 0 & A_1 & I_r & -I_r & 0 & 0 &
\cdots & A_1 & I_r & -I_r & 0 & 0\\
A_2 & 0 & 0 & I_s & -I_s & 0 & 0 & 0 & 0 & 0 & \cdots & 0 & 0 & 0 & 0 & 0\\
0 & 0 & 0 & 0 & 0 & A_2 & 0 & 0 & I_s & \!\!-I_s & \cdots & 0 & 0 & 0 & 0 & 0\\
\vdots & \vdots & \vdots & \vdots & \vdots & \vdots & \vdots & \vdots & \vdots
& \vdots & \ddots & \vdots & \vdots & \vdots & \vdots & \vdots\\
0 & 0 & 0 & 0 & 0 & 0 & 0 & 0 & 0 & 0 & \cdots & A_2 & 0 & 0 & I_s & -I_s\\
\end{array}
\right){\hat x}=b.
$$
This is again an $n$-fold integer program, with an $(r+s)\times (t+2r+2s)$
matrix ${\hat A}$, where ${\hat A_1}=(A_1,I_r,-I_r,0,0)$ and
${\hat A_2}=(A_2,0,0,I_s,-I_s)$.
Since $A$ is fixed, so is $\hat A$. It is now easy to extend the vector
$x\in\Z^{nt}$ determined above to a feasible point ${\hat x}$ of the
auxiliary program. Indeed, put ${\hat b}:=b-A^{(n)}x\in\Z^{r+ns}$;
now, for $i=1,\dots,r+ns$, simply choose an auxiliary variable ${\hat x_j}$
appearing only in the $i$-th equation, whose coefficient equals the sign
$\sign({\hat b}_i)$ of the corresponding entry of $\hat b$, and set
${\hat x_j}:=|{\hat b}_i|$. Define ${\hat w}\in\Z^{n(t+2r+2s)}$ by setting
${\hat w}:=0$ for each original variable and ${\hat w}:=-1$ for each
auxiliary variable, so that maximizing ${\hat w}{\hat x}$ is equivalent
to minimizing the sum of auxiliary variables. Now solve the auxiliary
linear integer program in polynomial time by applying the algorithm of Lemma \ref{NFoldAugmentation} corresponding to ${\hat A}$
to the data $n$, $\hat l$, $\hat u$, $\hat x$, and $\hat w$.
Since the auxiliary objective ${\hat w}{\hat x}$ is bounded above by zero,
the algorithm will output an optimal solution ${\hat x}^*$. If the optimal
objective value is negative, then the original $n$-fold program is
infeasible, whereas if the optimal value is zero, then the
restriction of ${\hat x}^*$ to the original variables is
a feasible point $x^*$ of the original integer program.
\epr

Combining Lemmas \ref{NFoldAugmentation} and \ref{Feasible}
we get at once the main result of this section.

\bt{NFoldTheorem}
For every fixed $(r+s)\times t$ integer matrix $A$ there is a polynomial
time algorithm that, given $n$, lower and upper bounds $l,u\in\Z_{\infty}^{nt}$, $w\in\Z^{nt}$,
and $b\in\Z^{r+ns}$, encoded as $[\l l,u,w,b\r]$, solves the
following linear $n$-fold integer programming problem,
$$\max\,\{wx\ :\ x\in\Z^{nt},\ A^{(n)}x=b,\ l\leq x\leq u\}\ .$$
\et

Again, as a special case of Theorem \ref{NFoldTheorem} we recover the
following result of \cite{DHOW} concerning linear integer programming
in standard form over $n$-fold matrices.

\bt{SpecialNFoldTheorem}
For every fixed $(r+s)\times t$ integer matrix $A$ there is a polynomial
time algorithm that, given $n$, linear functional $w\in\Z^{nt}$, and
right-hand side $b\in\Z^{r+ns}$, encoded as $[\l w,b\r]$,
solves the following linear $n$-fold integer program in standard form,
$$\max\,\{wx\ :\ x\in\N^{nt},\ A^{(n)}x=b\}\ .$$
\et

\subsection{Some Applications}
\label{salip}

\subsubsection{Three-Way Line-Sum Transportation Problems}
\label{twlstp}

Transportation problems form a very important class of discrete optimization
problems studied extensively in the operations research and mathematical
programming literature, see e.g. \cite{BR,KW,KLS,QS,Vla,YKK} and the
references therein. We will discuss this class of problem and its applications
to secure statistical data disclosure in more detail in \S \ref{mtpapisd}.

It is well known that $2$-way transportation problems
are polynomial time solvable, since they can be encoded
as linear integer programs over totally unimodular systems.
However, already $3$-way transportation problem are much more
complicated. Consider the following $3$-way transportation problem
over $p\times q \times n$ tables with all line-sums fixed,
$$\max\{wx\ :\ x\in\N^{p\times q\times n}\,,\ \sum_i x_{i,j,k}=z_{j,k}
\,,\ \sum_j x_{i,j,k}=v_{i,k}\,,\ \sum_k x_{i,j,k}=u_{i,j}\,\}\ .$$
The data for the problem consist of given integer numbers (lines-sums)
$u_{i,j}$, $v_{i,k}$, $z_{j,k}$ for $i=1,\dots,p$, $j=1,\dots,q$, $k=1,\dots,n$,
and a linear functional given by a $p\times q\times n$ integer array $w$
representing the transportation profit per unit on each cell. The problem is
to find a transportation, that is, a $p\times q\times n$ nonnegative integer
table $x$ satisfying the line sum constraints, which attains maximum profit
$wx=\sum_{i=1}^p\sum_{j=1}^q\sum_{k=1}^n w_{i,j,k}x_{i,j,k}$.

When at least two of the table sides, say $p,q$, are variable part
of the input, and even when the third side is fixed and as small as $n=3$,
this problem is already {\em universal} for integer programming
in a very strong sense \cite{DO2,DO4}, and in particular is NP-hard
\cite{DO1}; this will be discussed in detail and proved
in \S \ref{mtpapisd}. We now show that in contrast, when two sides,
say $p,q$, are fixed (but arbitrary), and one side $n$ is variable,
then the $3$-way transportation problem over such {\em long} tables
is an $n$-fold integer programming problem and therefore, as a consequence
of Theorem \ref{SpecialNFoldTheorem}, can be solved is polynomial time.

\bc{LinearThreeWay}
For every fixed $p$ and $q$ there is a polynomial time algorithm that,
given $n$, integer profit array $w\in\Z^{p\times q\times n}$,
and line-sums $u\in\Z^{p\times q}$, $v\in\Z^{p\times n}$ and
$z\in\Z^{q\times n}$, encoded as $[\l w,u,v,z \r]$,
solves the integer 3-way line-sum transportation problem
$$\max\{wx \ :\ x\in\N^{p\times q\times n}\,,\ \sum_i x_{i,j,k}=z_{j,k}
\,,\ \sum_j x_{i,j,k}=v_{i,k}\,,\ \sum_k x_{i,j,k}=u_{i,j}\,\}\ .$$
\ec
\bpr
Re-index $p\times q\times n$ arrays as $x=(x^1,\dots,x^n)$ with each component
indexed as $x^k:=(x^k_{i,j}):=(x_{1,1,k},\dots,x_{p,q,k})$ suitably indexed as
a $pq$ vector representing the $k$-th layer of $x$. Put $r:=t:=pq$ and $s:=p+q$,
and let $A$ be the $(r+s)\times t$ matrix with $A_1:=I_{pq}$ the $pq\times pq$
identity and with $A_2$ the $(p+q)\times pq$ matrix of equations of the
usual $2$-way transportation problem for $p\times q$ arrays. Re-arrange
the given line-sums in a vector $b:=(b^0,b^1,\dots,b^n)\in\Z^{r+ns}$
with $b^0:=(u_{i,j})$ and $b^k:=((v_{i,k}),(z_{j,k}))$ for $k=1,\dots,n$.

This translates the given $3$-way transportation problem into an
$n$-fold integer programming problem in standard form,
$$\max\,\{wx\ :\ x\in\N^{nt},\ A^{(n)}x=b\}\ ,$$
where the equations $A_1(\sum_{k=1}^n x^k)=b^0$ represent
the constraints $\sum_k x_{i,j,k}=u_{i,j}$ of all line-sums
where summation over layers occurs, and the equations $A_2x^k=b^k$
for $k=1,\dots,n$ represent the constraints $\sum_i x_{i,j,k}=z_{j,k}$
and $\sum_j x_{i,j,k}=v_{i,k}$ of all line-sums where summations are
within a single layer at a time.

Using the algorithm of Theorem \ref{SpecialNFoldTheorem},
this $n$-fold integer program, and hence the given $3$-way
transportation problem, can be solved in polynomial time.
\epr

\be{ExampleThreeWay}
We demonstrate the encoding of the $p\times q\times n$ transportation
problem as an $n$-fold integer program as in the proof of
Corollary \ref{LinearThreeWay} for $p=q=3$
(smallest case where the problem is genuinely $3$-dimensional).
Here we put $r:=t:=9$, $s:=6$, write
$$x^k:=(x_{1,1,k},x_{1,2,k},x_{1,3,k},x_{2,1,k}, x_{2,2,k},x_{2,3,k},
x_{3,1,k},x_{3,2,k},x_{3,3,k})\,,\ k=1,\dots,n\,,$$
and let the $(9+6)\times 9$ matrix
$A$ consist of $A_1=I_9$ the $9\times 9$ identity matrix and
$$A_2\quad:=\quad\left(
\begin{array}{ccccccccc}
  1 & 1 & 1 & 0 & 0 & 0 & 0 & 0 & 0 \\
  0 & 0 & 0 & 1 & 1 & 1 & 0 & 0 & 0 \\
  0 & 0 & 0 & 0 & 0 & 0 & 1 & 1 & 1 \\
  1 & 0 & 0 & 1 & 0 & 0 & 1 & 0 & 0 \\
  0 & 1 & 0 & 0 & 1 & 0 & 0 & 1 & 0 \\
  0 & 0 & 1 & 0 & 0 & 1 & 0 & 0 & 1 \\
\end{array}
\right)\quad.
$$
Then the corresponding $n$-fold integer program encodes
the $3\times 3\times n$ transportation problem as desired.
Already for this case, of $3\times 3\times n$ tables,
the only known polynomial time algorithm for the transportation
problem is the one underlying Corollary \ref{LinearThreeWay}.
\ee

Corollary \ref{LinearThreeWay} has a very broad generalization to
multiway transportation problems over long $k$-way tables of any
dimension $k$; this will be discussed in detail in \S \ref{mtpapisd}.

\subsubsection{Packing Problems and Cutting-Stock}
\label{ppacs}

We consider the following rather general class of packing problems which
concern maximum utility packing of many items of several types in various
bins subject to weight constraints. More precisely, the data is as follows.
There are $t$ types of items. Each item of type $j$ has integer weight $v_j$.
There are $n_j$ items of type $j$ to be packed. There are $n$ bins.
The weight capacity of bin $k$ is an integer $u_k$. Finally, there is a utility
matrix $w\in\Z^{t\times n}$ where $w_{j,k}$ is the utility of packing one item
of type $j$ in bin $k$. The problem is to find a feasible packing of maximum
total utility. By incrementing the number $t$ of types by $1$ and suitably
augmenting the data, we may assume that the last type $t$ represents
``slack items" which occupy the unused capacity in each bin, where the
weight of each slack item is $1$, the utility of packing any slack item
in any bin is $0$, and the number of slack items is the total residual
weight capacity $n_t:=\sum_{k=1}^n u_k-\sum_{j=1}^{t-1} n_jv_j$.
Let $x\in\N^{t\times n}$ be a variable matrix where $x_{j,k}$
represents the number of items of type $j$ to be packed in bin $k$.
Then the packing problem becomes the following linear integer program,
$$\max\{wx\ : x\in\N^{t\times n}
\,,\ \sum_j v_jx_{j,k}=u_k\,,\ \sum_k x_{j,k}=n_j\,\}\ .$$
We now show that this is in fact an $n$-fold integer programming problem and
therefore, as a consequence of Theorem \ref{SpecialNFoldTheorem}, can be solved
is polynomial time. While the number $t$ of types and type weights $v_j$ are
fixed, which is natural in many bin packing applications, the numbers $n_j$
of items of each type and the bin capacities $u_k$ may be very large.

\bc{Packing}
For every fixed number $t$ of types and integer type weights $v_1,\dots,v_t$,
there is a polynomial time algorithm that, given $n$ bins, integer
item numbers $n_1,\dots,n_t$,\break integer bin capacities $u_1,\dots,u_n$,
and $t\times n$ integer utility matrix $w$, encoded as\break
$[\l n_1,\dots,n_t,u_1,\dots,u_n,w\r]$,
solves the following integer bin packing problem,
$$\max\{wx\ : x\in\N^{t\times n}
\,,\ \sum_j v_jx_{j,k}=u_k\,,\ \sum_k x_{j,k}=n_j\,\}\ .$$
\ec
\bpr
Re-index the variable matrix as $x=(x^1,\dots,x^n)$ with
$x^k:=(x^k_1,\dots,x^k_t)$ where $x^k_j$ represents the number
of items of type $j$ to be packed in bin $k$ for al $j$ and $k$.
Let $A$ be the $(t+1)\times t$ matrix with $A_1:=I_t$ the $t\times t$
identity and with $A_2:=(v_1,\dots,v_t)$ a single row.
Re-arrange the given item numbers and bin capacities in a vector
$b:=(b^0,b^1,\dots,b^n)\in\Z^{t+n}$ with $b^0:=(n_1,\dots,n_t)$
and $b^k:=u_k$ for all $k$. This translates the bin packing problem into an
$n$-fold integer programming problem in standard form,
$$\max\,\{wx\ :\ x\in\N^{nt},\ A^{(n)}x=b\}\ ,$$
where the equations $A_1(\sum_{k=1}^n x^k)=b^0$ represent the
constraints $\sum_k x_{j,k}=n_j$ assuring that all items of each
type are packed, and the equations $A_2x^k=b^k$ for $k=1,\dots,n$
represent the constraints $\sum_j v_j x_{j,k}=u_k$ assuring
that the weight capacity of each bin is not exceeded
(in fact, the slack items make sure each bin is perfectly packed).

Using the algorithm of Theorem \ref{SpecialNFoldTheorem},
this $n$-fold integer program, and hence the given integer
bin packing problem, can be solved in polynomial time.
\epr

\be{CuttingStock}{\bf (cutting-stock problem).}
This is a classical manufacturing
problem \cite{GG}, where the usual setup is as follows:
a manufacturer produces rolls of material (such as scotch-tape or band-aid)
in one of $t$ different widths $v_1,\dots,v_t$.
The rolls are cut out from standard rolls of common large width $u$.
Given orders by customers for $n_j$ rolls of width $v_j$,
the problem facing the manufacturer is to meet the orders using the smallest
possible number of standard rolls. This can be cast as a bin packing problem
as follows. Rolls of width $v_j$ become items of type $j$ to be packed.
Standard rolls become identical bins, of capacity $u_k:=u$ each,
where the number of bins is set to be
$n:=\sum_{j=1}^t\lceil n_j/\lfloor u/v_j\rfloor\rceil$ which is
sufficient to accommodate all orders.
The utility of each roll of width $v_j$ is set to be its width negated
$w_{j,k}:=-v_j$ regardless of the standard roll $k$ from which it is cut
(paying for the width it takes). Introduce a new roll width $v_0:=1$,
where rolls of that width represent ``slack rolls" which occupy the
unused width of each standard roll, with utility $w_{0,k}:=-1$ regardless
of the standard roll $k$ from which it is cut (paying for the unused
width it represents), with the number of slack rolls set to be the total
residual width $n_0:=nu-\sum_{j=1}^t n_jv_j$. Then the cutting-stock problem
becomes a bin packing problem and therefore, by Corollary \ref{Packing},
for every fixed $t$ and fixed roll widths $v_1,\dots,v_t$, it is solvable
in time polynomial in $\sum_{j=1}^t\lceil n_j/\lfloor u/v_j\rfloor\rceil$
and $\l n_1,\dots,n_t,u \r$.

One common approach to the cutting-stock problem uses so-called
{\em cutting patterns}, which are feasible solutions of the knapsack
problem $\{y\in\N^t\,:\,\sum_{j=1}^t v_jy_j\leq u\}$. This is useful when
the common width $u$ of the standard rolls is of the same order of magnitude
as the demand role widths $v_j$. However, when $u$ is much larger than
the $v_j$, the number of cutting patterns becomes prohibitively large
to handle. But then the values $\lfloor u/v_j\rfloor$ are large and
hence $n:=\sum_{j=1}^t\lceil n_j / \lfloor u/v_j \rfloor \rceil$
is small, in which case the solution through the algorithm of
Corollary \ref{Packing} becomes particularly appealing.
\ee

\newpage
\section{Convex Integer Programming}
\label{cip}

In this section we discuss convex integer programming. In particular,
we extend the theory of \S \ref{lnfip} and show that convex $n$-fold
integer programming is polynomial time solvable as well.
In \S \ref{cipotus} we discuss convex integer programming over totally
unimodular matrices. In \S \ref{gbacip} we show the applicability of Graver
bases to convex integer programming. In \S \ref{cnfipipt} we combine
Theorem \ref{ConvexToLinear}, the results of \S \ref{lnfip}, and the
preparatory facts from \S \S \ref{gbacip}, and prove the main result of this
section, asserting that convex $n$-fold integer programming is polynomial
time solvable. We conclude with some applications in \S \ref{sacip}.

As in \S \ref{lnfip}, the feasible set $S$ is presented as the set of
integer points satisfying an explicitly given system of linear inequalities,
given in one of the forms
$$S:=\{x\in\N^n: Ax=b\}\quad \mbox{or}\quad
S:=\{x\in\Z^n:Ax=b,\ l\leq x\leq u\}\ \ ,$$
with matrix $A\in\Z^{m\times n}$, right-hand side $b\in\Z^m$, and
lower and upper bounds $l,u\in\Z_{\infty}^n$.

As demonstrated in \S \ref{l}, if the polyhedron
$P:=\{x\in\R^n:Ax=b,\ l\leq x\leq u\}$
is unbounded then the convex integer programming
problem with an oracle presented convex functional is rather hopeless.
Therefore, an algorithm that solves the convex integer programming problem
should either return an optimal solution, or assert that the program is
infeasible, or assert that the underlying polyhedron is unbounded.

Nonetheless, we do allow the lower and upper bounds $l,u$ to lie
in $\Z_{\infty}^n$ rather than $\Z^n$, since often the polyhedron
is bounded even though the variables are not bounded explicitly
(for instance, if each variable is bounded below only, and appears
in some equation all of whose coefficients are positive). This results
in broader formulation flexibility. Furthermore, in the next
subsections we prove auxiliary lemmas asserting that certain sets
cover all edge-directions of relevant polyhedra, which do hold also
in the unbounded case. So we now extend the notion of edge-directions,
defined in \S \ref{edaz} for polytopes, to polyhedra.
A {\em direction} of an edge ($1$-dimensional face) $e$ of a polyhedron $P$
is any nonzero scalar multiple of $y-x$ where $x,y$ are any two distinct
points in $e$. As before, a set {\em covers all edge-directions of $P$}
if it contains a direction of each edge of $P$.

\subsection{Convex Integer Programming over Totally Unimodular Systems}
\label{cipotus}

A matrix $A$ is {\em totally unimodular} if the determinant of every
square submatrix of $A$ lies in $\{-1,0,1\}$. Such matrices arise
naturally in network flows, ordinary ($2$-way) transportation problems,
and many other situations. A fundamental result in integer programming
\cite{HK} asserts that polyhedra defined by totally
unimodular matrices are integer. More precisely, if $A$ is an $m\times n$
totally unimodular matrix, $l,u\in\Z^n_{\infty}$, and $b\in\Z^m$, then
$$P_I:=\conv\{x\in\Z^n:Ax=b,\ l\leq x\leq u\}
\ =\ \{x\in\R^n:Ax=b,\ l\leq x\leq u\}:= P\ ,$$
that is, the underlying polyhedron $P$ coincides with its integer hull $P_I$.
This has two consequences useful in facilitating the solution of the
corresponding convex integer programming problem via the algorithm of Theorem \ref{ConvexToLinear}. First, the corresponding
linear integer programming problem can be solved by linear programming
over $P$ in polynomial time. Second, a set covering all
edge-directions of the implicitly given integer hull $P_I$, which is
typically very hard to determine, is obtained here as a set covering all
edge-directions of $P$ which is explicitly given and hence easier to determine.

We now describe a well known property of polyhedra of the above form.
A {\em circuit} of a matrix  $A\in\Z^{m\times n}$ is a nonzero primitive
minimal support element of $\L(A)$. So a circuit is a nonzero $c\in\Z^n$
satisfying $Ac=0$, whose entries are relatively prime integers, such that no
nonzero $c'$ with $Ac'=0$ has support strictly contained in the support of $c$.

\bl{Circuits}
For every $A\in\Z^{m\times n}$, $l,u\in\Z_{\infty}^n$, and $b\in\Z^m$,
the set of circuits of $A$ covers all edge-directions
of the polyhedron $P:=\{x\in\R^n:Ax=b,\ l\leq x\leq u\}$.
\el
\bpr
Consider any edge $e$ of $P$. Pick two distinct points $x,y\in e$
and set $g:=y-x$. Then $Ag=0$ and therefore, as can be easily proved by
induction on $|\supp(g)|$, there is a finite decomposition
$g=\sum_i \alpha_i c_i$ with $\alpha_i$ positive real number
and $c_i$ circuit of $A$ such that $\alpha_i c_i\sqsubseteq g$ for all $i$,
where $\sqsubseteq$ is the natural extension from $\Z^n$ to $\R^n$
of the partial order defined in \S \ref{gbalip}.
We claim that $x+\alpha_i c_i\in P$ for all $i$. Indeed, $c_i$ being a
circuit implies $A(x+\alpha_i c_i)=Ax=b$; and $l\leq x,x+g \leq u$ and
$\alpha_i c_i\sqsubseteq g$ imply $l\leq x+\alpha_i c_i\leq u$.

Now let $w\in\R^n$ be a linear functional uniquely maximized over
$P$ at the edge $e$. Then $w\alpha_i c_i=w(x+\alpha_i c_i)-wx\leq 0$
for all $i$. But $\sum (w\alpha_i c_i)=wg=wy-wx=0$, implying that in fact
$w\alpha_i c_i=0$ and hence $x+\alpha_i c_i\in e$ for all $i$. This implies
that each $c_i$ is a direction of $e$ (in fact, all $c_i$ are the same
and $g$ is a multiple of some circuit).
\epr

Combining Theorem \ref{ConvexToLinear} and Lemma \ref{Circuits}
we obtain the following statement.

\bt{TUMConvexIP}
For every fixed $d$ there is a polynomial time algorithm that, given
$m\times n$ totally unimodular matrix $A$, set $C\subset\Z^n$ containing
all circuits of $A$, vectors $l,u\in\Z_{\infty}^n$, $b\in\Z^m$, and $w_1,\dots,w_d\in\Z^n$,
and convex $c:\R^d\longrightarrow\R$ presented by a comparison oracle,
encoded as $[\l A,C,l,u,b,w_1,\dots,w_d\r]$, solves the convex integer program
$$\max\, \{c(w_1x,\dots, w_dx)\ :\ x\in\Z^n,\ Ax=b,\ l\leq x\leq u\}\ .$$
\et
\bpr
First, check in polynomial time using linear programming whether the
objective function of any of the following $2n$ linear programs is unbounded,
$$\max\, \{\pm y_i:y\in P\},\ i=1,\dots,n,
\quad P:=\{y\in\R^n\ :\ Ay=b,\ l\leq y\leq u\}\ .$$
If any is unbounded then terminate, asserting that $P$ is unbounded.
Otherwise, let $\rho$ be the least integer upper bound on the absolute value
of all optimal objective values. Then $P\subseteq [-\rho,\rho]^n$ and
$S:=\{y\in\Z^n:Ay=b,\ l\leq y\leq u\}\subset P$ is finite of radius
$\rho(S)\leq \rho$. In fact, since $A$ is totally unimodular,
$P_I=P=\conv(S)$ and hence $\rho(S)=\rho$. Moreover, by Cramer's rule,
$\l\rho\r$ is polynomially bounded in $\l A,l,u,x \r$.

Now, since $A$ is totally unimodular, using linear programming over $P_I=P$ we
can simulate in polynomial time a linear discrete optimization oracle for $S$.
By Lemma \ref{Circuits}, the given set $C$, which contains all circuits of $A$,
also covers all edge-directions of $\conv(S)=P_I=P$. Therefore we can apply
the algorithm of Theorem \ref{ConvexToLinear} and solve the given
convex $n$-fold integer programming problem in polynomial time.
\epr

While the number of circuits of an $m\times n$ matrix $A$ can be as large as
$2{n\choose m+1}$ and hence exponential in general, it is nonetheless
relatively small in that it is bounded in terms of $m$ and $n$ only and
is independent of the matrix $A$ itself. Furthermore, it may happen that
the number of circuits is much smaller than the upper bound $2{n\choose m+1}$.
Also, if in a class of matrices, $m$ grows slowly in terms of $n$, say
$m=O(\log n)$, then this bound is subexponential. In such situations,
the above theorem may provide a good strategy for solving
convex integer programming over totally unimodular systems.

\subsection{Graver Bases and Convex Integer Programming}
\label{gbacip}

We now extend the statements of \S \ref{cipotus} about totally unimodular
matrices to arbitrary integer matrices. The next lemma shows that the
Graver basis of any integer matrix covers all edge-directions of the
integer hulls of polyhedra defined by that matrix.

\bl{GraverEdges}
For every $A\in\Z^{m\times n}$, $l,u\in\Z_{\infty}^n$, and $b\in\Z^m$,
the Graver basis $\G(A)$ of $A$ covers all edge-directions
of the polyhedron $P_I:=\conv\{x\in\Z^n:Ax=b,\ l\leq x\leq u\}$.
\el
\bpr
Consider any edge $e$ of $P_I$ and pick two distinct points $x,y\in e\cap\Z^n$.
Then $g:=y-x$ is in $\L(A)\setminus\{0\}$. Therefore, by Lemma \ref{Conformal},
there is a conformal sum $g=\sum_i h_i$ with $h_i\in\G(A)$ for all $i$.
We claim that $x+h_i\in P_I$ for all $i$. Indeed, first note
that $h_i\in\G(A)\subset\L(A)$ implies $Ah_i=0$ and hence $A(x+h_i)=Ax=b$;
and second note that $l\leq x,x+g \leq u$ and $h_i\sqsubseteq g$ imply
that $l\leq x+h_i\leq u$.

Now let $w\in\Z^n$ be a linear functional uniquely maximized over
$P_I$ at the edge $e$. Then $wh_i=w(x+h_i)-wx\leq 0$ for all $i$.
But $\sum (wh_i)=wg=wy-wx=0$, implying that in fact $wh_i=0$
and hence $x+h_i\in e$ for all $i$. Therefore each $h_i$
is a direction of $e$ (in fact, all $h_i$ are the same
and $g$ is a multiple of some Graver basis element).
\epr

Combining Theorems \ref{ConvexToLinear} and \ref{GraverIP}
and Lemma \ref{GraverEdges} we obtain the following statement.

\bt{GraverConvexIP}
For every fixed $d$ there is a polynomial time algorithm that, given integer
$m\times n$ matrix $A$, its Graver basis $\G(A)$, $l,u\in\Z_{\infty}^n$,
$x\in\Z^n$ with $l\leq x\leq u$, $w_1,\dots,w_d\in\Z^n$, and convex
$c:\R^d\longrightarrow\R$ presented by a comparison oracle,
encoded as $[\l A,\G(A),l,u,x,w_1,\dots,w_d\r]$,
solves the convex integer program with $b:=Ax$,
$$\max\, \{c(w_1z,\dots, w_dz)\ :\ z\in\Z^n,\ Az=b,\ l\leq z\leq u\}\ .$$
\et
\bpr
First, check in polynomial time using linear programming whether the
objective function of any of the following $2n$ linear programs is unbounded,
$$\max\, \{\pm y_i:y\in P\},\ i=1,\dots,n,
\quad P:=\{y\in\R^n\ :\ Ay=b,\ l\leq y\leq u\}\ .$$
If any is unbounded then terminate, asserting that $P$ is unbounded.
Otherwise, let $\rho$ be the least integer upper bound on the absolute value
of all optimal objective values. Then $P\subseteq [-\rho,\rho]^n$ and
$S:=\{y\in\Z^n:Ay=b,\ l\leq y\leq u\}\subset P$ is finite of radius
$\rho(S)\leq \rho$. Moreover, by Cramer's rule, $\l\rho\r$
is polynomially bounded in $\l A,l,u,x \r$.

Using the given Graver basis and applying the algorithm of
Theorem \ref{GraverIP} we can simulate in polynomial time a linear
discrete optimization oracle for $S$. Furthermore, by Lemma \ref{GraverEdges},
the given Graver basis covers all edge-directions of the integer hull
$P_I:=\conv\{y\in\Z^n:Ay=b, l\leq y\leq u\}=\conv(S)$. Therefore we can apply
the algorithm of Theorem \ref{ConvexToLinear} and solve the
given convex program in polynomial time.
\epr

\subsection{Convex N-fold Integer Programming in Polynomial Time}
\label{cnfipipt}

We now extend the result of Theorem \ref{NFoldTheorem} and show that convex
integer programming problems over $n$-fold systems can be solved in polynomial
time as well. As explained in the beginning of this section, the algorithm
either returns an optimal solution, or asserts that the program is
infeasible, or asserts that the underlying polyhedron is unbounded.

\bt{NFoldConvex}
For every fixed $d$ and fixed $(r+s)\times t$ integer matrix $A$ there
is a polynomial time algorithm that, given $n$, lower and upper bounds
$l,u\in\Z_{\infty}^{nt}$, $w_1,\dots,w_d\in\Z^{nt}$,\break $b\in\Z^{r+ns}$,
and convex functional $c:\R^d\longrightarrow\R$ presented
by a comparison oracle, encoded as $[\l l,u,w_1,\dots,w_d,b \r]$,
solves the convex $n$-fold integer programming problem
$$\max\,\{c(w_1x,\dots,w_dx)\ :\ x\in\Z^{nt},\ A^{(n)}x=b,\ l\leq x\leq u\}\ .$$
\et
\bpr
First, check in polynomial time using linear programming whether the
objective function of any of the following $2nt$ linear programs is unbounded,
$$\max\, \{\pm y_i:y\in P\},\ i=1,\dots,nt,
\quad P:=\{y\in\R^{nt}\ :\ A^{(n)}y=b,\ l\leq y\leq u\}\ .$$
If any is unbounded then terminate, asserting that $P$ is unbounded.
Otherwise, let $\rho$ be the least integer upper bound on the absolute value
of all optimal objective values. Then $P\subseteq [-\rho,\rho]^{nt}$ and
$S:=\{y\in\Z^{nt}:A^{(n)}y=b,\ l\leq y\leq u\}\subset P$ is finite of radius
$\rho(S)\leq \rho$. Moreover, by Cramer's rule, $\l\rho\r$
is polynomially bounded in $n$ and $\l l,u,b \r$.

Using the algorithm of Theorem \ref{NFoldTheorem} we can simulate
in polynomial time a linear discrete optimization oracle for $S$.
Also, using the algorithm of Theorem \ref{GraverComputation}
we can compute in polynomial time the Graver basis $\G(A^{(n)})$ which,
by Lemma \ref{GraverEdges}, covers all edge-directions of
$P_I:=\conv\{y\in\Z^{nt}:A^{(n)}y=b,l\leq y\leq u\}=\conv(S)$.
Therefore we can apply the algorithm of Theorem \ref{ConvexToLinear}
and solve the given convex $n$-fold integer programming problem
in polynomial time.
\epr

Again, as a special case of Theorem \ref{NFoldConvex} we recover the
following result of \cite{DHORW} concerning convex integer programming
in standard form over $n$-fold matrices.

\bt{SpecialNFoldConvex}
For every fixed $d$ and fixed $(r+s)\times t$ integer matrix $A$ there is a
polynomial time algorithm that, given $n$, linear functionals
$w_1,\dots,w_d\in\Z^{nt}$, right-hand side\break $b\in\Z^{r+ns}$, and convex
functional $c:\R^d\longrightarrow\R$ presented by a comparison oracle,
encoded as $[\l w_1,\dots,w_d,b \r]$, solves the
convex $n$-fold integer program in standard form
$$\max\,\{c(w_1x,\dots,w_dx)\ :\ x\in\N^{nt},\ A^{(n)}x=b\}\ .$$
\et

\subsection{Some Applications}
\label{sacip}

\subsubsection{Transportation Problems and Packing Problems}
\label{tpapp}

Theorems \ref{NFoldConvex} and \ref{SpecialNFoldConvex} generalize
Theorems \ref{NFoldTheorem} and \ref{SpecialNFoldTheorem} by broadly
extending the class of objective functions that can be maximized in
polynomial time over $n$-fold systems. Therefore all applications
discussed in \S \ref{salip} automatically extend accordingly.

First, we have the following analog of Corollary \ref{LinearThreeWay}
for the {\em convex integer transportation problem} over long $3$-way tables.
This has a very broad further generalization to multiway transportation
problems over long $k$-way tables of any dimension $k$, see \S \ref{mtpapisd}.

\bc{ConvexThreeWay}
For every fixed $d,p,q$ there is a polynomial time algorithm that,
given $n$, arrays $w_1,\dots,w_d\in\Z^{p\times q\times n}$,
line-sums $u\in\Z^{p\times q}$, $v\in\Z^{p\times n}$ and
$z\in\Z^{q\times n}$, and convex functional $c:\R^d\longrightarrow\R$
presented by a comparison oracle, encoded as $[\l w_1,\dots,w_d,u,v,z \r]$,
solves the convex integer 3-way line-sum transportation problem
\begin{eqnarray*}
\max\{\,c(w_1x,\dots,w_dx) \ &:& \ x\in\N^{p\times q\times n}\,,\\
& & \sum_i x_{i,j,k}=z_{j,k}\,,\ \sum_j x_{i,j,k}=v_{i,k}
\,,\ \sum_k x_{i,j,k}=u_{i,j}\,\}\ .
\end{eqnarray*}
\ec

Second, we have the following analog of
Corollary \ref{Packing} for convex bin packing.

\bc{ConvexPacking}
For every fixed $d$, number of types $t$, and type weights $v_1,\dots,v_t\in\Z$,
there is a polynomial time algorithm that, given $n$ bins, item
numbers $n_1,\dots,n_t\in\Z$, bin capacities $u_1,\dots,u_n\in\Z$,
utility matrices $w_1,\dots,w_d\in\Z^{t\times n}$, and convex
functional $c:\R^d\longrightarrow\R$ presented by a comparison oracle,
encoded as $[\l n_1,\dots,n_t,u_1,\dots,u_n,w_1,\dots,w_d\r]$,
solves the convex integer bin packing problem,
$$\max\{\,c(w_1x,\dots,w_dx)\ : x\in\N^{t\times n}
\,,\ \sum_j v_jx_{j,k}=u_k\,,\ \sum_k x_{j,k}=n_j\,\}\ .$$
\ec

\subsubsection{Vector Partitioning and Clustering}
\label{vpac}

The vector partition problem concerns the partitioning of $n$
items among $p$ players to maximize social value subject to
constraints on the number of items each player can receive.
More precisely, the data is as follows. With each item $i$ is
associated a vector $v_i\in\Z^k$ representing its utility under $k$ criteria.
The utility of player $h$ under ordered partition $\pi=(\pi_1,\dots,\pi_p)$ of the
set of items $\{1,\dots,n\}$ is the sum $v^{\pi}_h:=\sum_{i\in\pi_h} v_i$
of utility vectors of items assigned to $h$ under $\pi$.
The social value of $\pi$ is the balancing
$c(v^{\pi}_{1,1},\dots,v^{\pi}_{1,k},\dots,v^{\pi}_{p,1},\dots,v^{\pi}_{p,k})$
of the player utilities, where $c$ is a convex functional on $\R^{pk}$.
In the constrained version, the partition must be of a given {\em shape},
i.e. the number $|\pi_h|$ of items that player $h$ gets is required to be a
given number $\lambda_h$ (with $\sum\lambda_h=n$). In the unconstrained
version, there is no restriction on the number of items per player.

Vector partition problems have applications in diverse areas such
as load balancing, circuit layout, ranking, cluster analysis,
inventory, and reliability, see e.g. \cite{BHR,BH,FOR,HOR,OSc}
and the references therein. Here is a typical example.

\be{Clustering}{\bf (minimal variance clustering).}
This problem has numerous applications in the analysis of statistical
data: given $n$ observed points $v_1,\dots,v_n$ in $k$-space,
group them into $p$ clusters $\pi_1,\dots,\pi_p$
that minimize the sum of cluster variances given by
$$\sum_{h=1}^p{1\over|\pi_h|}\sum_{i\in \pi_h}
||v_i-({1\over|\pi_h|}\sum_{i\in\pi_h} v_i)||^2\ .$$
Consider instances where there are $n=pm$ points and the
desired clustering is balanced, that is, the clusters should
have equal size $m$. Suitable manipulation of the sum of variances
expression above shows that the problem is equivalent to a constrained
vector partition problem, where $\lambda_h=m$ for all $h$, and where the
convex functional $c:\R^{pk}\longrightarrow\R$
(to be maximized) is the Euclidean norm squared, given by
$$c(z)\ =\ ||z||^2\ =\ \sum_{h=1}^p\sum_{i=1}^k|z_{h,i}|^2\ .$$
\ee

If either the number of criteria $k$ or the number of players $p$ is variable,
the partition problem is intractable since it instantly captures NP-hard
problems \cite{HOR}. When both $k,p$ are fixed, both the constrained and
unconstrained versions of the vector partition problem are polynomial time
solvable \cite{HOR,OSc}. We now show that vector partition problems
(either constrained or unconstrained) are in fact convex $n$-fold integer
programming problems and therefore, as a consequence of
Theorem \ref{SpecialNFoldConvex}, can be solved is polynomial time.

\bc{Partition}
For every fixed number $p$ of players and number $k$ of criteria,
there is a polynomial time algorithm that, given $n$, item vectors
$v_1,\dots,v_n\in\Z^k$, $\lambda_1,\dots,\lambda_p\in\N$,\break and convex
functional $c:\R^{pk}\longrightarrow\R$ presented by a comparison oracle,
encoded as $[\l v_1,\dots,v_n,\lambda_1,\dots,\lambda_p\r]$, solves
the constrained and unconstrained partitioning problems.
\ec
\bpr
There is an obvious one-to-one correspondence between partitions and
matrices $x\in\{0,1\}^{p\times n}$ with all column-sums equal to one,
where partition $\pi$ corresponds to the
matrix $x$ with $x_{h,i}=1$ if $i\in\pi_h$ and $x_{h,i}=0$ otherwise.
Let $d:=pk$ and define $d$ matrices $w_{h,j}\in\Z^{p\times n}$
by setting $(w_{h,j})_{h,i}:=v_{i,j}$ for all $h=1,\dots,p$,
$i=1,\dots,n$ and $j=1,\dots,k$, and setting all other entries to zero.
Then for any partition $\pi$ and its corresponding matrix $x$ we have
$v^{\pi}_{h,j}=w_{h,j}x$ for all $h=1,\dots,p$ and $j=1,\dots,k$.
Therefore, the unconstrained vector partition
problem is the convex integer program
$$\max\{\,c(w_{1,1} x, \dots, w_{p,k} x) \ :\ x\in\N^{p\times n}
\,,\ \sum_h x_{h,i}=1\,\}\ .$$
Suitably arranging the variables in a vector, this becomes
a convex $n$-fold integer program with a $(0+1)\times p$
defining matrix $A$, where $A_1$ is empty and $A_2:=(1,\dots,1)$.

Similarly, the constrained vector partition
problem is the convex integer program
$$\max\{\,c(w_{1,1} x, \dots, w_{p,k} x) \ :\ x\in\N^{p\times n}
\,,\ \sum_h x_{h,i}=1\,,\ \sum_i x_{h,i}=\lambda_h\,\}\ .$$
This again is a convex $n$-fold integer program, now with a
$(p+1)\times p$ defining matrix $A$, where now $A_1:=I_p$ is the
$p\times p$ identity matrix and $A_2:=(1,\dots,1)$ as before.

Using the algorithm of Theorem \ref{SpecialNFoldConvex},
this convex $n$-fold integer program, and hence the given vector
partition problem, can be solved in polynomial time.
\epr

\newpage
\section{Multiway Transportation Problems and Privacy in Statistical Databases}
\label{mtpapisd}

Transportation problems form a very important class of discrete
optimization problems. The feasible points in a transportation problem
are the multiway tables (``contingency tables" in statistics) such that
the sums of entries over some of their lower dimensional sub-tables such as
lines or planes (``margins" in statistics) are specified. Transportation
problems and their corresponding transportation polytopes have been used and
studied extensively in the operations research and mathematical programming
literature, as well as in the statistics literature in the context of
secure statistical data disclosure and management by public agencies,
see \cite{AT,BR,Cox,DT,DLTZ,KW,KLS,QS,Vla,YKK} and references therein.

In this section we completely settle the algorithmic complexity of treating
multiway tables and discuss the applications to transportation problems
and secure statistical data disclosure, as follows.
After introducing some terminology in \S \ref{tam}, we go on to describe,
in \S \ref{tut}, a universality result that shows that
``short" $3$-way $r\times c\times 3$ tables, with variable number $r$
of rows and variable number $c$ of columns but fixed small number $3$ of
layers (hence ``short"), are {\em universal} in a very strong sense.
In \S \ref{tcotmtp} we discuss the general multiway transportation problem.
Using the results of \S \ref{tut} and the results on linear and convex
$n$-fold integer programming from \S \ref{lnfip} and \S \ref{cip},
we show that the transportation problem is intractable for short
$3$-way $r\times c\times 3$ tables but polynomial time treatable
for ``long" $(k+1)$-way $m_1\times\cdots\times m_k\times n$ tables,
with $k$ and the sides $m_1,\ldots, m_k$ fixed (but arbitrary),
and the number $n$ of layers variable (hence ``long").
In \S \ref{paeu} we turn to discuss data privacy and security and
consider the central problem of detecting entry uniqueness in
tables with disclosed margins. We show that as a consequence
of the results of \S \ref{tut} and \S \ref{tcotmtp}, and in analogy
to the complexity of the transportation problem established
in \S \ref{tcotmtp}, the entry uniqueness problem is intractable for
short $3$-way $r\times c\times 3$ tables but polynomial time decidable
for long $(k+1)$-way $m_1\times\cdots\times m_k\times n$ tables.

\subsection{Tables and Margins}
\label{tam}

We start with some terminology on tables,
margins and transportation polytopes.\break
A {\em $k$-way table} is an $m_1\times\cdots\times m_k$
array $x=(x_{i_1,\dots,i_k})$ of nonnegative integers.
A {\em $k$-way transportation polytope} (or simply {\em $k$-way polytope}
for brevity) is the set of all $m_1\times\cdots\times m_k$
nonnegative arrays $x=(x_{i_1,\dots,i_k})$ such that the sums
of the entries over some of their lower dimensional sub-arrays
(margins) are specified. More precisely, for any tuple $(i_1,\dots,i_k)$
with $i_j\in\{1,\dots,m_j\}\cup\{+\}$, the corresponding {\em margin}
$x_{i_1,\dots,i_k}$ is the sum of entries of $x$ over all coordinates
$j$ with $i_j=+$. The {\em support} of $(i_1,\dots,i_k)$ and of
$x_{i_1,\dots,i_k}$ is the set $\supp(i_1,\dots,i_k):=\{j:i_j\neq +\}$
of non-summed coordinates. For instance, if $x$ is a $4\times5\times3\times2$
array then it has $12$ margins with support $F=\{1,3\}$ such as
$x_{3,+,2,+}=\sum_{i_2=1}^5\sum_{i_4=1}^2 x_{3,i_2,2,i_4}$.
A collection of margins is {\em hierarchical} if, for some family $\F$ of
subsets of $\{1,\dots,k\}$, it consists of all margins $u_{i_1,\dots,i_k}$
with support in $\F$. In particular, for any $0\leq h\leq k$,
the collection of {\em all $h$-margins} of $k$-tables is the hierarchical
collection with $\F$ the family of all $h$-subsets of $\{1,\dots,k\}$.
Given a hierarchical collection of margins $u_{i_1,\dots,i_k}$
supported on a family $\F$ of subsets of $\{1,\dots,k\}$,
the corresponding {\em $k$-way polytope}
is the set of nonnegative arrays with these margins,
$$T_{\F} \ :=\  \left\{\,x\in\R_+^{m_1\times\cdots \times m_k}
\ : \ x_{i_1,\dots,i_k}\,=\,u_{i_1,\dots, i_k}
\,,\ \ \supp(i_1,\dots,i_k)\in\F\,\right\}\ .$$
The integer points in this polytope are precisely the $k$-way
tables with the given margins.

\subsection{The Universality Theorem}
\label{tut}

We now describe the following {\em universality} result of \cite{DO2,DO4}
which shows that, quite remarkably, {\em any} rational polytope is a
short $3$-way $r\times c\times 3$ polytope with all line-sums specified.
(In the terminology of \S \ref{tam} this is the $r\times c\times 3$ polytope
$T_{\F}$ of all $2$-margins fixed, supported
on the family $\F=\{\{1,2\},\{1,3\},\{2,3\}\}$.)
By saying that a polytope $P\subset\R^p$ is {\em representable} as a polytope
$Q\subset\R^q$ we mean in the strong sense that there is an injection
$\sigma:\{1,\dots,p\}\longrightarrow\{1,\dots,q\}$
such that the coordinate-erasing projection
$$\pi:\R^q\longrightarrow\R^p \ :\ x=(x_1,\dots,x_q)\mapsto
\pi(x)=(x_{\sigma(1)},\dots,x_{\sigma(p)})$$ provides a bijection
between $Q$ and $P$ and between the sets of integer points $Q\cap\Z^q$
and $P\cap\Z^p$. In particular, if $P$ is representable as $Q$ then
$P$ and $Q$ are isomorphic in any reasonable sense: they are linearly
equivalent and hence all linear programming related problems over the two are
polynomial time equivalent; they are combinatorially equivalent and hence
they have the same face numbers and facial structure; and they are integer
equivalent and therefore all integer programming and integer counting
related problems over the two are polynomial time equivalent as well.

We provide only an outline of the proof of the following statement; complete
details and more consequences of this theorem can be found in \cite{DO2,DO4}.

\bt{Universality}
There is a polynomial time algorithm that, given
$A\in\Z^{m\times n}$ and $b\in\Z^m$, encoded as $[\l A,b\r]$,
produces $r,c$ and line-sums $u\in\Z^{r\times c}$,
$v\in\Z^{r\times 3}$ and $z\in\Z^{c\times 3}$ such that the polytope
$P:=\{y\in\R_+^n:Ay=b\}$ is representable as the $3$-way polytope
$$T\quad:=\quad\{x\in\R_+^{r\times c\times 3}\ :\ \sum_i x_{i,j,k}=z_{j,k}
\,,\ \sum_j x_{i,j,k}=v_{i,k}\,,\ \sum_k x_{i,j,k}=u_{i,j}\,\}\ .$$
\et
\bpr
The construction proving the theorem consists of three
polynomial time steps, each representing a polytope of
a given format as a polytope of another given format.

First, we show that any $P:=\{y\geq 0:Ay=b\}$ with $A,b$ integer
can be represented
in polynomial time as $Q:=\{x\geq 0:Cx=d\}$ with $C$ matrix all entries
of which are in $\{-1,0,1,2\}$. This reduction of coefficients will enable
the rest of the steps to run in polynomial time.
For each variable $y_j$ let
$k_j:=\max\{\lfloor \log_2 |a_{i,j}| \rfloor : i=1,\dots m\}$ be the maximum
number of bits in the binary representation of the absolute value of any entry
$a_{i,j}$ of $A$. Introduce variables $x_{j,0},\dots,x_{j,k_j}$, and relate
them by the equations $2x_{j,i}-x_{j,i+1}=0$. The representing injection
$\sigma$ is defined by $\sigma(j):=(j,0)$, embedding $y_j$ as $x_{j,0}$.
Consider any term $a_{i,j}\,y_j$ of the original system.
Using the binary expansion $|a_{i,j}|=\sum_{s=0}^{k_j}t_s 2^s$ with all
$t_s\in\{0,1\}$, we rewrite this term as $\pm\sum_{s=0}^{k_j}t_s x_{j,s}$.
It is not hard to verify that this represents $P$ as $Q$ with
defining $\{-1,0,1,2\}$-matrix.

Second, we show that any $Q:=\{y\geq 0:Ay=b\}$ with $A,b$ integer
can be represented as a face $F$ of a $3$-way polytope with all plane-sums
fixed, that is, a face of a $3$-way polytope $T_{\F}$ of all $1$-margins
fixed, supported on the family $\F=\{\{1\},\{2\},\{3\}\}$.

Since $Q$ is a polytope and hence bounded, we can compute (using Cramer's rule)
an integer upper bound $U$ on the value of any coordinate $y_j$ of any
$y\in Q$. Note also that a face of a $3$-way polytope $T_{\F}$ is the set
of all $x=(x_{i,j,k})$ with some entries forced to zero; these entries are
termed ``forbidden", and the other entries are termed ``enabled".

For each variable $y_j$, let $r_j$ be the largest between the sum of
positive coefficients of $y_j$ and the sum of absolute values of
negative coefficients of $y_j$ over all equations,
$$r_j\quad:=\quad \max
\left(\sum_k\{a_{k,j}:a_{k,j}>0\}\,,\,\sum_k\{|a_{k,j}|:a_{k,j}<0\}\right)\ .$$
Assume that $A$ is of size $m\times n$.
Let $r:=\sum_{j=1}^n r_j$, $R:=\{1,\dots,r\}$, $h:=m+1$ and $H:=\{1,\dots,h\}$.
We now describe how to construct vectors $u,v\in\Z^r, z\in\Z^h$,
and a set $E\subset R\times R\times H$ of triples - the enabled,
non-forbidden, entries - such that the polytope $Q$ is represented as the
face $F$ of the corresponding $3$-way polytope of $r\times r\times h$ arrays
with plane-sums $u,v,z$ and only entries indexed by $E$ enabled,
\begin{eqnarray*}
F \ :=\ \{ \,x\in\R_+^{r\times r\times h}\ & : &
\ x_{i,j,k}=0\ \, \mbox{for all}\ \, (i,j,k)\notin E\,,\,
\ \, \mbox{and}\ \, \\
& &  \sum_{i,j} x_{i,j,k}=z_k\,,
\ \sum_{i,k} x_{i,j,k}=v_j\,,\ \sum_{j,k} x_{i,j,k}=u_i\, \}\ .
\end{eqnarray*}
We also indicate the injection
$\sigma:\{1,\dots,n\}\longrightarrow R\times R\times H$
giving the desired embedding of coordinates $y_j$ as coordinates
$x_{i,j,k}$ and the representation of $Q$ as $F$.

Roughly, each equation $k=1,\dots,m$ is encoded in a ``horizontal plane"
$R\times R\times\{k\}$\break (the last plane $R\times R\times\{h\}$ is included
for consistency with its entries being ``slacks"); and
each variable $y_j$, $j=1,\dots,n$ is encoded in a ``vertical box"
$R_j\times R_j\times H$, where $R=\biguplus_{j=1}^n R_j$ is the natural
partition of $R$ with $|R_j|=r_j$ for all $j=1,\dots,n$, that is,
with $R_j:=\{1+\sum_{l<j}r_l,\dots,\sum_{l\leq j}r_l\}$.

Now, all ``vertical" plane-sums are set to the same value $U$, that is,
$u_j:=v_j:=U$ for $j=1,\dots,r$. All entries not in the union
$\biguplus_{j=1}^n R_j\times R_j\times H$
of the variable boxes will be forbidden. We now describe the enabled entries
in the boxes; for simplicity we discuss the box $R_1\times R_1\times H$,
the others being similar. We distinguish between the two cases
$r_1=1$ and $r_1\geq 2$. In the first case, $R_1=\{1\}$; the box, which
is just the single line $\{1\}\times\{1\}\times H$, will have exactly two
enabled entries $(1,1,k^+),(1,1,k^-)$ for suitable\break $k^+$, $k^-$ to be
defined later. We set $\sigma(1):=(1,1,k^+)$, namely embed $y_1=x_{1,1,k^+}$.
We define the {\em complement} of the variable $y_1$ to be
${\bar y}_1:=U-y_1$ (and likewise for the other variables). The vertical
sums $u,v$ then force ${\bar y}_1=U-y_1=U-x_{1,1,k^+}=x_{1,1,k^-}$, so the
complement of $y_1$ is also embedded. Next, consider the case $r_1\geq 2$.
For each $s=1,\dots, r_1$, the line $\{s\}\times \{s\}\times H$
(respectively, $\{s\}\times \{1+ (s \mod r_1)\}\times H$) will contain
one enabled entry $(s,s,k^+(s))$ (respectively, $(s,1+ (s \mod r_1),k^-(s))$.
All other entries of $R_1\times R_1\times H$ will be forbidden.
Again, we set $\sigma(1):=(1,1,k^+(1))$, namely embed $y_1=x_{1,1,k^+(1)}$;
it is then not hard to see that, again, the vertical sums $u,v$ force
$x_{s,s,k^+(s)}=x_{1,1,k^+(1)}=y_1$ and
$x_{s,1+ (s \small \mod r_1),k^-(s)}=U-x_{1,1,k^+(1)}={\bar y}_1$
for each $s=1,\dots, r_1$. Therefore, both $y_1$ and ${\bar y}_1$
are each embedded in $r_1$ distinct entries.

We now encode the equations by defining the horizontal
plane-sums $z$ and the indices $k^+(s), k^-(s)$ above as follows.
For $k=1,\dots,m$, consider the $k$-th equation\break $\sum_j a_{k,j}y_j=b_k$.
Define the index sets $J^+:=\{j:a_{k,j}>0\}$ and $J^-:=\{j:a_{k,j}<0\}$,
and set $z_k:=b_k+U\cdot\sum_{j\in J^-}|a_{k,j}|$.
The last coordinate of $z$ is set for consistency
with $u,v$ to be $z_h=z_{m+1}:=r\cdot U-\sum_{k=1}^m z_k$.
Now, with ${\bar y}_j:=U-y_j$ the complement of
variable $y_j$ as above, the $k$-th equation can be rewritten as
$$\sum_{j\in J^+} a_{k,j}y_j+\sum_{j\in J^-} |a_{k,j}|{\bar y_j}
=\sum_{j=1}^n a_{k,j}y_j+U\cdot \sum_{j\in J^-}|a_{k,j}|
=b_k+U\cdot\sum_{j\in J^-}|a_{k,j}|=z_k.$$
To encode this equation, we simply ``pull down" to the
corresponding $k$-th horizontal plane as many copies of each variable
$y_j$ or ${\bar y}_j$ by suitably setting $k^+(s):=k$ or $k^-(s):=k$.
By the choice of $r_j$ there are sufficiently
many, possibly with a few redundant copies which are absorbed in
the last hyperplane by setting $k^+(s):=m+1$ or $k^-(s):=m+1$.
This completes the encoding and provides the desired representation.

Third, we show that any $3$-way polytope with plane-sums fixed
and entry bounds,
$$F := \{\,y\in\R_+^{l\times m\times n}  : \sum_{i,j} y_{i,j,k}=c_k\,,
\sum_{i,k} y_{i,j,k}=b_j\,, \sum_{j,k} y_{i,j,k}=a_i\,,
\ y_{i,j,k}\leq e_{i,j,k} \}\, ,$$
can be represented as a $3$-way polytope with line-sums fixed
(and no entry bounds),
$$T\ :=\ \{\,x\in\R_+^{r\times c\times 3}\ :\ \sum_I
x_{I,J,K}=z_{J,K}\,,\ \sum_J x_{I,J,K}=v_{I,K}\,,\ \sum_K
x_{I,J,K}=u_{I,J}\,\}\ .$$
In particular, this implies that any face $F$ of a $3$-way polytope
with plane-sums fixed can be represented as a $3$-way polytope $T$ with
line-sums fixed: forbidden entries are encoded by setting a
``forbidding" upper-bound $e_{i,j,k}:=0$ on all forbidden entries
$(i,j,k)\notin E$ and an ``enabling" upper-bound $e_{i,j,k}:=U$
on all enabled entries $(i,j,k)\in E$. We describe the presentation,
but omit the proof that it is indeed valid; further details on this step can
be found in \cite{DO1,DO2,DO4}. We give explicit formulas
for $u_{I,J}, v_{I,K}, z_{J,K}$ in terms of $a_i, b_j, c_k$ and $e_{i,j,k}$
as follows. Put $r:=l\cdot m$ and $c:=n+l+m$. The first index $I$ of each
entry $x_{I,J,K}$ will be a pair $I=(i,j)$ in the $r$-set
$$\{(1,1),\dots,(1,m),(2,1),\dots,(2,m),\dots,(l,1),\dots,(l,m)\}\ .$$
The second index $J$ of each entry $x_{I,J,K}$ will be a pair $J=(s,t)$
in the $c$-set
$$\{(1,1),\dots,(1,n),(2,1),\dots,(2,l),(3,1),\dots,(3,m)\}\ .$$ The
last index $K$ will simply range in the $3$-set $\{1,2,3\}$.  We
represent $F$ as $T$ via the injection $\sigma$ given explicitly by
$\sigma(i,j,k):=((i,j),(1,k),1)$, embedding each variable $y_{i,j,k}$
as the entry $x_{(i,j),(1,k),1}$. Let $U$ now denote the minimal
between the two values $\max\{a_1,\dots,a_l\}$ and
$\max\{b_1,\dots,b_m\}$.  The line-sums ($2$-margins) are set to be
$$u_{(i,j),(1,t)}= e_{i,j,t},\quad u_{(i,j),(2,t)}=\left\{ \begin{array}{ll}
U & \mbox {if $t=i$,} \\
   0 & \mbox{otherwise.} \cr
\end{array} \right. , \quad u_{(i,j),(3,t)}=\left\{ \begin{array}{ll} U &
\mbox{if $t=j,$} \\ 0 & \mbox{otherwise.} \cr \end{array} \right.$$

$$v_{(i,j),t}=\left\{ \begin{array}{ll} U & \mbox{if $t=1$,} \\
                                       e_{i,j,+} & \mbox{if $t=2$,}\\
                                       U & \mbox{if $t=3$.} \cr
\end{array} \right. ,\quad
z_{(i,j),1}=\left\{ \begin{array}{ll} c_j & \mbox{if $i=1$,} \\
                                      m\cdot U-a_j & \mbox{if $i=2$,}\\
                                       0 & \mbox{if $i=3$.} \cr
\end{array} \right.
$$
$$z_{(i,j),2}=\left\{ \begin{array}{ll} e_{+,+,j}-c_j & \mbox{if $i=1$,} \\
                                      0 & \mbox{if $i=2$,}\\
                                       b_j & \mbox{if $i=3$.} \cr
\end{array} \right. , \quad
z_{(i,j),3}=\left\{ \begin{array}{ll} 0 & \mbox{if $i=1$,} \\
                                      a_j & \mbox{if $i=2$,}\\
                                       l\cdot U-b_j & \mbox{if $i=3$.} \cr
\end{array} \right.\quad.
$$

Applying the first step to the given rational polytope $P$, applying
the second step to the resulting $Q$, and applying the third step to the
resulting $F$, we get in polynomial time a $3$-way $r\times c\times 3$
polytope $T$ of all line-sums fixed representing $P$ as claimed.
\epr

\subsection{The Complexity of the Multiway Transportation Problem}
\label{tcotmtp}

We are now finally in position to settle the complexity of the general
multiway transportation problem. The data for the problem
consists of: positive integers $k$ (table dimension) and $m_1,\dots,m_k$
(table sides); family $\F$  of subsets of $\{1,\dots,k\}$
(supporting the hierarchical collection of margins to
be fixed); integer values $u_{i_1,\dots,i_k}$ for all margins supported
on $\F$; and integer ``profit" $m_1\times\cdots\times m_k$ array $w$.
The transportation problem is to find an $m_1\times\cdots\times m_k$
table having the given margins and attaining maximum profit, or assert than
none exists. Equivalently, it is the linear integer programming problem
of maximizing the linear functional defined by $w$ over the
transportation polytope $T_{\F}$,
$$\max\left\{wx\,:\,x\in\N^{m_1\times\cdots \times m_k}
\ : \ x_{i_1,\dots,i_k}\,=\,u_{i_1,\dots, i_k}
\,,\ \ \supp(i_1,\dots,i_k)\in\F\,\right\}\ .$$

The following result of \cite{DO1} is an immediate consequence
of Theorem \ref{Universality}. It asserts that if two sides of
the table are variable part of the input then the transportation
problem is intractable already for short $3$-way tables with
$\F=\{\{1,2\},\{1,3\},\{2,3\}\}$ supporting all $2$-margins (line-sums).
This result can be easily extended to $k$-way tables of any dimension
$k\geq 3$ and $\F$ the collection of all $h$-subsets of $\{1,\dots,k\}$
for any $1<h<k$ as long as two sides of the table are variable;
we omit the proof of this extended result.

\bc{NPCThreeWay}
It is NP-complete to decide, given $r,c$, and line-sums
$u\in\Z^{r\times c}$,\break $v\in\Z^{r\times 3}$, and $z\in\Z^{c\times 3}$,
encoded as $[\l u,v,z \r]$, if the following set of tables is nonempty,
$$S\ :=\ \{x\in\N^{r\times c\times 3}\ :\ \sum_i x_{i,j,k}=z_{j,k}
\,,\ \sum_j x_{i,j,k}=v_{i,k}\,,\ \sum_k x_{i,j,k}=u_{i,j}\,\}\ .$$
\ec
\bpr
The integer programming feasibility problem is to decide, given
$A\in\Z^{m\times n}$ and $b\in\Z^m$, if $\{y\in\N^n:Ay=b\}$ is nonempty.
Given such $A$ and $b$, the polynomial time algorithm of
Theorem \ref{Universality} produces $r,c$ and $u\in\Z^{r\times c}$,
$v\in\Z^{r\times 3}$, and $z\in\Z^{c\times 3}$,
such that $\{y\in\N^n:Ay=b\}$ is nonempty if and only if the set $S$
above is nonempty. This reduces integer programming feasibility to short
$3$-way line-sum transportation feasibility. Since the former is
NP-complete (see e.g. \cite{Sch}), so turns out to be the latter.
\epr

We now show that in contrast, when all sides but one are fixed
(but arbitrary), and one side $n$ is variable, then the corresponding
long $k$-way transportation problem for any hierarchical collection
of margins is an $n$-fold integer programming problem and therefore,
as a consequence of Theorem \ref{SpecialNFoldTheorem}, can be solved
is polynomial time. This extends Corollary \ref{LinearThreeWay}
established in \S \ref{twlstp} for $3$-way line-sum transportation.

\bc{LinearKWay}
For every fixed $k$, table sides $m_1,\dots,m_k$, and family
$\F$ of subsets of $\{1,\dots,k+1\}$, there is a polynomial time
algorithm that, given $n$, integer values\break
$u=(u_{i_1,\dots,i_{k+1}})$ for all margins supported on $\F$,
and integer $m_1\times\cdots\times m_k\times n$ array $w$,
encoded as $[\l u,w\r]$, solves the linear
integer multiway transportation problem
$$\max\left\{wx\,:\,x\in\N^{m_1\times\cdots\times m_k\times n}
,\ x_{i_1,\dots,i_{k+1}}\,=\,u_{i_1,\dots, i_{k+1}}
,\ \supp(i_1,\dots,i_{k+1})\in\F\,\right\}\,.$$
\ec
\bpr
Re-index the arrays as $x=(x^1,\dots,x^n)$ with each
$x^j=(x_{i_1,\dots,i_k,j})$ a suitably indexed $m_1m_2\cdots m_k$ vector
representing the $j$-th layer of $x$. Then the transportation problem
can be encoded as an $n$-fold integer programming problem in standard form,
$$\max\,\{wx\ :\ x\in\N^{nt},\ A^{(n)}x=b\}\ ,$$
with an $(r+s)\times t$ defining matrix $A$ where $t:=m_1m_2\cdots m_k$
and $r,s$, $A_1$ and $A_2$ are determined from $\F$,
and with right-hand side $b:=(b^0,b^1,\dots,b^n)\in\Z^{r+ns}$
determined from the margins $u=(u_{i_1,\dots,i_{k+1}})$,
in such a way that the equations $A_1(\sum_{j=1}^n x^j)=b^0$ represent
the constraints of all margins $x_{i_1,\dots,i_k,+}$
(where summation over layers occurs), whereas the equations $A_2x^j=b^j$
for $j=1,\dots,n$ represent the constraints of all margins
$x_{i_1,\dots,i_k,j}$ with $j\neq +$
(where summations are within a single layer at a time).

Using the algorithm of Theorem \ref{SpecialNFoldTheorem},
this $n$-fold integer program, and hence the given multiway
transportation problem, can be solved in polynomial time.
\epr

The proof of Corollary \ref{LinearKWay} shows that the set of feasible points
of any long $k$-way transportation problem, with all sides but one fixed and
one side $n$ variable, for any hierarchical collection of margins,
is an $n$-fold integer programming problem. Therefore, as a consequence of
Theorem \ref{SpecialNFoldConvex}, we also have the following extension
of Corollary \ref{LinearKWay} for the convex integer multiway
transportation problem over long $k$-way tables.

\bc{ConvexKWay}
For every fixed $d$, $k$, table sides $m_1,\dots,m_k$, and family $\F$ of
subsets of $\{1,\dots,k+1\}$, there is a polynomial time algorithm that,
given $n$, integer values $u=(u_{i_1,\dots,i_{k+1}})$ for all margins
supported on $\F$, integer $m_1\times\cdots\times m_k\times n$ arrays
$w_1,\dots,w_d$, and convex functional $c:\R^d\longrightarrow\R$
presented by a comparison oracle, encoded as $[\l u,w_1,\dots,w_d\r]$,
solves the convex integer multiway transportation problem
\begin{eqnarray*}
\max\,\{\, c(w_1x,\dots,w_dx)\ &:&
\ x\in\N^{m_1\times\cdots\times m_k\times n}\,,\\
& & x_{i_1,\dots,i_{k+1}}\ =\ u_{i_1,\dots, i_{k+1}}
\,,\ \supp(i_1,\dots,i_{k+1})\in\F\,\}\ .
\end{eqnarray*}
\ec

\subsection{Privacy and Entry-Uniqueness}
\label{paeu}

A common practice in the disclosure of a multiway table containing
sensitive data is to release some of the table margins rather than the
table itself, see e.g. \cite{Cox,DT,DLTZ} and the references therein.
Once the margins are released, the security of any specific entry of
the table is related to the set of possible values that can occur in
that entry in any table having the same margins as those
of the source table in the data base. In particular, if this
set consists of a unique value, that of the source table, then this
entry can be exposed and privacy can be violated. This raises the following
fundamental {\em entry-uniqueness problem}: given a consistent disclosed
(hierarchical) collection of margin values, and a specific entry index,
is the value that can occur in that entry in any table having these margins
unique? We now describe the results of \cite{Onn2} that
settle the complexity of this problem, and interpret the
consequences for secure statistical data disclosure.

First, we show that if two sides of the table are variable part of the
input then the entry-uniqueness problem is intractable already for
short $3$-way tables with all $2$-margins (line-sums) disclosed
(corresponding to $\F=\{\{1,2\},\{1,3\},\{2,3\}\}$). This can be easily
extended to $k$-way tables of any dimension $k\geq 3$ and $\F$ the
collection of all $h$-subsets of $\{1,\dots,k\}$ for any $1<h<k$ as
long as two sides of the table are variable; we omit the proof of this
extended result. While this result indicates that the disclosing agency
may not be able to check for uniqueness, in this situation, some
consolation is in that an adversary will be computationally unable
to identify and retrieve a unique entry either.

\bc{HardUniqueness}
It is coNP-complete to decide, given $r,c$, and line-sums
$u\in\Z^{r\times c}$,\break $v\in\Z^{r\times 3}$, $z\in\Z^{c\times 3}$, encoded
as $[\l u,v,z \r]$, if the entry $x_{1,1,1}$ is the same in all tables in\break
$$\{x\in\N^{r\times c\times 3}\ :\ \sum_i x_{i,j,k}=z_{j,k}
\,,\ \sum_j x_{i,j,k}=v_{i,k}\,,\ \sum_k x_{i,j,k}=u_{i,j}\,\}\ .$$
\ec
\bpr
The {\em subset-sum problem}, well known to be NP-complete, is the following:
given positive integers $a_0,a_1,\dots,a_m$, decide if there is an
$I\subseteq\{1,\dots,m\}$ with $a_0=\sum_{i\in I}a_i$. We reduce the complement
of subset-sum to entry-uniqueness. Given $a_0,a_1,\dots,a_m$, consider
the polytope in $2(m+1)$ variables $y_0,y_1\dots,y_m,z_0,z_1,\dots,z_m$,
$$P\ :=\ \{(y,z)\in\R_+^{2(m+1)}\ :
\ a_0 y_0-\sum_{i=1}^m a_i y_i=0\,,\ y_i+z_i=1\,,\ i=0,1\dots,m\,\}\ .$$
First, note that it always has one integer point with $y_0=0$,
given by $y_i=0$ and $z_i=1$ for all $i$.
Second, note that it has an integer point with $y_0\neq 0$ if and only if
there is an $I\subseteq\{1,\dots,m\}$ with $a_0=\sum_{i\in I}a_i$,
given by $y_0=1$, $y_i=1$ for $i\in I$, $y_i=0$ for
$i\in\{1,\dots,m\}\setminus I$, and $z_i=1-y_i$ for all $i$.
Lifting $P$ to a suitable $r\times c\times 3$ line-sum polytope $T$
with the coordinate $y_0$ embedded in the entry $x_{1,1,1}$ using
Theorem \ref{Universality}, we find that $T$ has a table with $x_{1,1,1}=0$,
and this value is unique among the tables in $T$ if and only if there
is {\em no} solution to the subset-sum problem with $a_0,a_1,\dots,a_m$.
\epr

Next we show that, in contrast, when all table sides but one are fixed
(but arbitrary), and one side $n$ is variable, then, as
a consequence of Corollary \ref{LinearKWay}, the corresponding
long $k$-way entry-uniqueness problem for any hierarchical collection
of margins can be solved is polynomial time. In this situation,
the algorithm of Corollary \ref{EasyUniqueness} below allows disclosing
agencies to efficiently check possible collections of margins before
disclosure: if an entry value is not unique then disclosure may be
assumed secure, whereas if the value is unique then disclosure may
be risky and fewer margins should be released. Note that this situation,
of long multiway tables, where one category is significantly richer
than the others, that is, when each sample point can take many values
in one category and only few values in the other categories, occurs
often in practical applications, e.g., when one category is the individuals
age and the other categories are binary (``yes-no"). In such situations,
our polynomial time algorithm below allows disclosing agencies
to check entry-uniqueness and make learned decisions on secure disclosure.

\bc{EasyUniqueness}
For every fixed $k$, table sides $m_1,\dots,m_k$, and family $\F$ of subsets of
$\{1,\dots,k+1\}$, there is a polynomial time algorithm that, given $n$, integer
values\break $u=(u_{j_1,\dots,j_{k+1}})$ for all margins supported on $\F$,
and entry index $(i_1,\dots,i_{k+1})$,\break encoded as $[n,\l u\r]$, decides
if the entry $x_{i_1,\dots,i_{k+1}}$ is the same in all tables in the set
$$\{x\in\N^{m_1\times\cdots\times m_k\times n}
\ :\ x_{j_1,\dots,j_{k+1}}\ =\ u_{j_1,\dots, j_{k+1}}
\,,\ \supp(j_1,\dots,j_{k+1})\in\F\,\}\ .$$
\ec
\bpr
By Theorem \ref{LinearKWay} we can solve in polynomial time
both transportation problems
$$l\quad:=\quad\min\left\{x_{i_1,\dots,i_{k+1}}
\ :\ x\in\N^{m_1\times\cdots\times m_k\times n}\,,\quad x\in T_{\F}\right\}\ ,$$
$$u\quad:=\quad\max\left\{x_{i_1,\dots,i_{k+1}}
\ :\ x\in\N^{m_1\times\cdots\times m_k\times n}\,,\quad x\in T_{\F}\right\}\ ,$$
over the corresponding $k$-way transportation polytope
$$T_{\F}\,:=\,\left\{\,x\in\R_+^{m_1\times\cdots \times m_k\times n}
\,:\, x_{j_1,\dots,j_{k+1}}\,=\,u_{j_1,\dots,j_{k+1}}
\,,\ \ \supp(j_1,\dots,j_{k+1})\in\F\,\right\}\,.$$
Clearly, entry $x_{i_1,\dots,i_{k+1}}$ has the same value in all
tables with the given (disclosed) margins if and only if $l=u$,
completing the description of the algorithm and the proof.
\epr


\begin{thebibliography}{00}
\bibitem{KLS}
Kleinschmidt, P., Lee, C.W., Schannath, H.: Transportation problems which
can be solved by the use of Hirsch-paths for the dual problems.
Math. Prog. {\bf 37} (1987) 153--168
\end{thebibliography}
\end{document}